\documentclass{article}
\usepackage{amsmath,amsthm,amssymb,mathrsfs}
\usepackage{fancyhdr}

\usepackage{pdfsync}
\usepackage{amsmath}
\usepackage{amsfonts}
\usepackage{graphicx}
\usepackage[all,cmtip]{xy}
\usepackage{color}
\allowdisplaybreaks

\def\p#1#2{\frac{\partial #1}{\partial #2}}

\def\ds{\displaystyle}
\newfont{\Blackboard}{msbm10 scaled 1200}

\newfont{\roma}{cmr10 scaled 1200}

\providecommand{\bysame}{\leavevmode\hbox to3em{\hrulefill}\thinspace}
\providecommand{\MR}{\relax\ifhmode\unskip\space\fi MR }

\providecommand{\href}[2]{#2}

\def\a{\alpha}

\newtheorem{theorem}{{}\hskip\parindent Theorem}[section]
\newtheorem{lem}{{}\hskip\parindent Lemma}[section]
\newtheorem{proposition}{{}\hskip\parindent Proposition}[section]
\newtheorem{example}{{}\hskip\parindent Example}[section]

\newtheorem{definition}{{}\hskip\parindent Definition}[section]
\newtheorem{rem}{{}\hskip\parindent Remark}[section]

\def\p{\prime}

\def\beq{\arraycolsep=1.5pt\begin{eqnarray}}
\def\eeq{\end{eqnarray}}

\large
\title{Second-Order Necessary Optimality Conditions for Multi-Objective Optimal Control Problems  on Riemannian Manifolds \thanks{This work is  supported by the National Science
Foundation of China under grant 12371451, and the Natural Science Foundation of Sichuan province under grants 2025ZNSFSC0077. 
}
\date{}
\author{Li Deng\thanks{School of Mathematics,   Southwest Jiaotong University, Chengdu 611756, Sichuan Province, China. {\small\it E-mail:} {\small\tt
dengli@swjtu.edu.cn}.}}
}

\begin{document}
\maketitle

\begin{quote}
\begin{small}
{\bf Abstract} \,\,\,
In this paper, we investigate the multi-objective optimal control problem of ordinary differential equations on Riemannian manifolds. We first obtain the second-order necessary conditions for weak Pareto optimal solutions for multi-objective optimal control problems with fixed terminal time, and then extend these results to multi-objective optimal control problems with free terminal time, deriving the corresponding second-order necessary conditions for weak Pareto optimal solutions. Our main results (i.e., Theorem \ref{sen} and Theorem \ref{snt}) show that weak Pareto optimal solutions depend on the curvature tensor of the Riemannian manifold. Finally, we provide an example (i.e., Example \ref{exg}) as an application of our main results, illustrating how our findings differ from existing related results (see Remark \ref{reer}).
\\[3mm]
{\bf Keywords}\,\,\, multi-objective optimal control problem, control systems on Riemannian manifolds, weak Pareto optimal solution, second-order necessary condition
\\[3mm]
{\bf MSC (2010) \,\,\,
49K15, 49K30, 70Q05, 90C29} \\[3mm]
\end{small}
\end{quote}

\setcounter{equation}{0}

\section{Introduction }
\def\theequation{1.\arabic{equation}}

 Let $n\in\mathbb N$ and $M$ be a complete, simply connected, $n$-dimensional manifold with a Riemannian metric $g$. Let $\nabla$ be the Levi-Civita connection on $M$ with respect to $g$, $\rho(\cdot,\cdot)$ be the distance function on $M$, $T_xM$ be the tangent space of  $M$ at $x\in M$, and $T^*_x M$ be the  cotangent space. Denote by $\langle\cdot,\cdot\rangle$ and $|\cdot|$ the inner product and the norm over $T_xM$ with respect to $g$, respectively.  We also denote by $|\cdot|$ the norm over $T_x^*M$, which is defined by $|\eta|\stackrel{\triangle}{=}\sup\{\eta(X)|\,X\in T_xM,\,|X|\leq 1\}$ (see \cite[(2.5)]{cdz}). As usual, denote the tangent and the cotangent bundles on $M$ respectively by
 $$ T M \equiv\{(x,X)|\; x\in M, X\in T_xM\} \ \ \hbox{ and } \ \ T^* M \equiv \{(x,\eta)|\;x\in M,
 \eta\in T_x^*M\}.$$
 Also, denote by $\mathcal X(M)$ and $C^\infty(M)$ the set of smooth vector fields and the set of smooth functions on $M$, respectively. For any given $h\in C^\infty(M)$,  $d h$  stands for the differential of $h$.

Let $j,m\in\mathbb N$, $k\in\mathbb N$, $r\in\mathbb N$,  $U\subseteq \mathbb R^m$, and the maps $f:[0,+\infty)\times M\times U\to T M$,  $\phi_0=(\phi_{0_1},\cdots,\phi_{0_r})^\top: M\times M\to\mathbb R^r$,  $\phi=(\phi_1,\cdots,\phi_j)^\top: M\times M \to \mathbb R^{j}$  and $\psi: M\times M \to \mathbb R^k$   satisfy suitable assumptions to be given later.

First, we consider an optimal control problem defined over a fixed time interval. Given $T>0$, the multi-objective  optimal control problem is stated as follows: 
\begin{description}
\item[$(MP)_{r,T}$]  Minimize $
J( u(\cdot), x(\cdot)) \equiv \phi_0(x(0), x(T))
$
subject to 
the state equation:
\beq\label{25}
\dot{x}(t)=f(t,x(t),u(t)),\quad x(t)\in M\quad a.e. \,t\in[0,T],
\eeq
the control constraint 
$
u(\cdot)\in\mathcal U_T\equiv \big\{u(\cdot):[0,T]\to U|\,u(\cdot) \mbox{ is measurable}\big\}$, 
and
endpoints-constraints:
\begin{equation}\label{n22}
\phi_i(x(0), x(T))\leq 0,\;i=1,\cdots,j\quad\textrm{and} \quad\psi(x(0), x(T))=O_{\mathbb R^k},
\end{equation}
where  $\dot{x}(t)\stackrel{\triangle}{=}\frac{d}{dt}x(t)
$ for $t\in [0,T]$, and $O_{\mathbb R^k}$ is the origin in $\mathbb R^k$. 
\end{description}

For  problem $(MP)_{r,T}$, we need to provide the definition of an optimal solution.

\begin{definition}\label{adoT}
For problem $(MP)_{r,T}$, $(u(\cdot), x(\cdot))$ is called an admissible solution, and $u(\cdot)$ is called an admissible control, if $(u(\cdot), x(\cdot))$ satisfies (\ref{25}), $u(\cdot)\in\mathcal U_T$, and (\ref{n22}). An admissible solution $(\bar u(\cdot),\bar x(\cdot))$ is referred to as a weak Pareto optimal solution,
if there exists no admissible solution $(u(\cdot), x(\cdot))$ such that $\phi_{0_i}(x(0),x(T)) < \phi_{0_i}(\bar x(0),\bar x(T))$ for all $ i=1,\cdots,r$.
In this case, $\bar u(\cdot)$ and  $\bar x(\cdot)$ are respectively called a weak Pareto optimal control and a weak Pareto optimal state.

\end{definition}


Then, we consider multi-objective time optimal control problem. 
\begin{description}
\item[$(MP)_r$]  Minimize $
\hat J( u(\cdot), x(\cdot), T) \equiv \phi_0(x(0), x(T))
$
subject to 
the state equation 
\beq\label{25t}
\dot{x}(t)=f(t,x(t),u(t)),\quad x(t)\in M\quad a.e. \,t\in[0,+\infty),
\eeq the  constraint $u(\cdot)\in\mathcal U\equiv \big\{u(\cdot):[0,+\infty)\to U|\; u(\cdot) \mbox{ is measurable}\big\}$, $T>0$
and endpoints-constraints (\ref{n22}).
\end{description}

Similar to problem $(MP)_{r,T}$, we also need to provide the definition of an optimal solution for problem $(MP)_r$.

\begin{definition}\label{ado}
For problem $(MP)_r$,
a triple $(u(\cdot), x(\cdot), T)$ is called an admissible solution, if it satisfies (\ref{n22}),  (\ref{25t}), $T>0$ and $u(\cdot)\in\mathcal U$. Here $u(\cdot)$ is called an admissible control. Moreover, an admissible solution $(\bar u(\cdot), \bar x(\cdot),\bar T)$ is said to be a weak Pareto optimal solution, if there exists no admissible solution $(u(\cdot), x(\cdot),T)$ such that $\phi_{0_i}(x(0),x(T))
<\phi_{0_i}(\bar x(0),\bar x(T))$ for $i=1,\cdots,r$.
 Here, $\bar u(\cdot)$, $\bar x(\cdot)$ and $\bar T$ are called a weak Pareto  optimal control, a weak Pareto optimal state and a weak Pareto optimal time, respectively.
\end{definition}

 The multi-objective optimal control problem is a generalization of the single-objective optimal control problem.
  When 
$r=1$, the problem $(MP)_{1,T}$ (or $(MP)_1$) becomes the traditional optimal control problem. In this case, the definition of a  weak Pareto optimal solution  is consistent with the definition of an optimal solution  in the traditional sense.

In this paper, we are concerned with the second-order necessary conditions for weak Pareto optimal solutions of problem $(MP)_{r,T}$ and problem $(MP)_r$. 

Why do we study control systems on differentiable manifolds? There are three reasons. First, many control systems based on classical mechanics theory are control systems on differentiable manifolds. For example, the monograph \cite{Bloch} lists some control problems of mechanical systems with motion constraints, which are precisely control systems on differentiable manifolds. Second, from the perspective of General Relativity, the four-dimensional spacetime in which we live is actually curved. Therefore, for the characterization of control systems on a large scale, the state should be valued in curved space. Finally, for control systems where the state is pointwise constrained to a set, studying their optimal control problems is often challenging. However, if this set happens to be a differentiable manifold (such as the control system (\ref{exs})), then the control system can be regarded as a control system on the differentiable manifold, thereby formally eliminating the pointwise state constraints.

The research results on the necessary optimality conditions for multi-objective  optimal control problems involving ordinary differential equations are primarily focused on control systems of ordinary differential equations in Euclidean spaces. Among them, the literature \cite{kwy, BA04, zhu00} discusses the first-order necessary conditions for Pareto optimal solutions, while the literature \cite{KQWY20, kty,  kb23} derives the second-order necessary conditions for weak Pareto optimal solutions. Specifically, \cite{KQWY20} deals with systems having free terminal time and where the control and state (with respect to time) satisfy pointwise inequality-type constraints,  \cite{kb23} addresses control systems where the state and control (with respect to time) satisfy pointwise inequality-type constraints,  and \cite{kty} investigates control systems where the control and state (with respect to time) satisfy pointwise inequality-type constraints.

Compared with existing results on the second-order necessary conditions for multi-objective optimal control problems, our study of the problems $(MP)_{r,T}$ and $(MP)_r$ presents two distinct characteristics. Firstly, the state in our control system is constrained pointwise in time on a differentiable manifold.  This type of constraint poses significant challenges when using variational tools to study optimality conditions.  Secondly, the state of our control system must satisfy equality-type constraints at the terminal time (i.e., the second constraint in (\ref{n22})). Such constraints often introduce difficulties in studying the second-order necessary optimality conditions.

To address the first challenge, we consider control systems on Riemannian manifolds. Using variational techniques to study the optimality conditions, we apply tools from Riemannian geometry to perform variations on the control systems. Specifically, we find that the second-order variation of the state in the control system depends on the curvature tensor of the Riemannian manifold, and consequently, the derived second-order necessary conditions also depend on the curvature tensor of the Riemannian manifold (see Theorem \ref{sen} and Theorem \ref{snt}). To overcome the second difficulty, we utilize  Brouwer's fixed-point theorem and the separation theorem of convex sets. Notably, in order to apply the Brouwer's fixed-point theorem, we need to construct a continuous mapping. The continuity of this mapping requires that the control constraint set (i.e., $U$) is convex.

Finally, as an application of our results, we provide an example (i.e., Example \ref{exg}). Through this example, we can see the differences between our results and those in the aforementioned literature, as detailed in Remark \ref{reer}.

This paper is organized as follows. In Section \ref{mr2}, we present the main results (namely, Theorem \ref{sen} and Theorem \ref{snt}) and provide an illustrative example (Example \ref{exg}). The proofs of the main results are detailed in Section \ref{pmr3}. In Section \ref{ae21}, we derive the local coordinate representations of some key expressions involved in Theorem \ref{sen} and Theorem \ref{snt}, and apply these results to Example \ref{exg}.

\setcounter{equation}{0}

\section{Main results }\label{mr2}
\def\theequation{2.\arabic{equation}}

\subsection{Notations and assumptions}

First, we shall introduce some notations that will be used throughout this paper.

Denote by $\mathcal T_s^r(x)$ for $x\in M$ and $r,s\in\mathbb N\cup\{0\}$ the tensor space of type $(r,s)$ at $x$, and by $\mathcal T_s^r(M)$ the set of all tensor fields of type $(r,s)$ over $M$.
We also denote   by $i(x)$, $\rho(x,y)$,  $| F(x)|$ and $\nabla F$ the injectivity radius at $x\in M$, the distance between $x$ and $y$, the norm of the tensor field $F$ at $x$ and the covariant differential of the tensor field $F$ of the Riemannian manifold $(M,g)$, respectively.
For each $(t,u)\in[0,T]\times U$, if   $\xi(t,\cdot,u)$ defines a tensor field on $M$ (see \cite[Section 2.2]{cdz}), then for each $x\in M$,  $\nabla_X\xi(t,\cdot,u)$  represents the covariant derivative     of tensor field $\xi(t,\cdot,u)$ with respect  to the tangent vector $X\in T_xM$ (see \cite[(2.9)]{cdz}). Additionally,  $\nabla_x\xi(t,\cdot,u)$ denotes the covariant differential of this tensor field at $x$ (see \cite[(2.10)]{cdz}).  In particular, when $\xi$ is a function, we denote by $d_x\xi(t,\cdot,u)$ the differential of $\xi(t,\cdot,u)$ (see \cite[p.42]{pth}). For a function $\zeta:M\times M\to\mathbb R$, we denote respectively by $\nabla_i\zeta(x_1,x_2)$ and $d_i\zeta(x_1,x_2)$ the covariant differential  and differential of $\zeta$ with respect to $x_i$, where $i=1,2$. 
Moreover, we define higher order covariant derivatives  of functions defined on $M\times M$. For $i,j\in\{1,2\}$ with $i \not= j$,  $(x_1,x_2)\in M\times M$ and $X,Y,Z\in \mathcal X( M)$, we set
\begin{align}\label{sodn}\begin{array}{l}
\nabla_i\nabla_j\zeta(x_1,x_2)(X,Y)\stackrel{\triangle}{=}Y(x_i)\Big(X(x_j)\big(\zeta(x_1,x_2)\big)\Big)=Y(x_i)\Big(\nabla_j\zeta(x_1,x_2)\big(X(x_j)\big)\Big),
\\[2mm] \nabla_i^2\zeta(x_1,x_2)(X,Y)=Y(x_i)\Big(X(x_i)\left(\zeta(x_1,x_2)\right)\Big)-\nabla_{Y(x_i)}X\big(\zeta(x_1,x_2)\big),
\end{array}
\end{align}
where we have used \cite[(5.9)]{dzgJ} and the definition of  the differential of a $\mathbb R$-valued function (see \cite[p.42]{pth}). 
We refer to \cite[Section 2]{cdz} and  the references cited therein for the notions mentioned above.

More generally,  for a vector-valued map $\xi=(\xi_1,\cdots,
 \xi_l)^\top: M\times M\to\mathbb R^l$ with $l\in\mathbb N$, the corresponding covariant   derivatives $\nabla_1\xi$ and $\nabla_2\xi$ are defined as follows:
\begin{align}\label{fdv}\begin{array}{r}
\nabla_i\xi(x_1,x_2)(X)\equiv \big(\nabla_i\xi_1(x_1,x_2)(X),\cdots,\nabla_i\xi_l(x_1,x_2)(X)\big)^\top,\quad\forall\,x_1,x_2\in M,
\\  X\in T_{x_i}M,\;i=1,2,
\end{array}\end{align}
 and the corresponding norm is given by
 \begin{align}\label{cdv}
|\nabla_i\xi(x_1,x_2)|=\sum_{j=1}^l|\nabla_i\xi_j(x_1,x_2)|,\quad i=1,2.
\end{align}

We also set
\begin{align}\label{sdv}\begin{array}{l}
\nabla_i^2\xi(x_1,x_2)(X,Y)=\big(\nabla_i^2\xi_1(x_1,x_2)(X,Y),\cdots,\nabla_i^2\xi_l(x_1,x_2)(X,Y)\big)^\top,
\\[2mm] \nabla_i\nabla_j\xi(x_1,x_2)(X,Y)=\big(\nabla_i\nabla_j\xi_1(x_1,x_2)(X,Y),\cdots,\nabla_i\nabla_j\xi_l(x_1,x_2)(X,Y)\big)^\top,
\end{array}
\end{align}
for all $X,Y\in\mathcal X(M)$. Particularly, for a weak Pareto optimal solution $(\bar u(\cdot),\bar x(\cdot))$, we set
\begin{align}\label{msd}\begin{array}{ll}
D^2\xi(\bar x(0),\bar x(T))(V)
=&\nabla_1^2\xi(\bar x(0),\bar x(T))(V(\bar x(0)), V(\bar x(0)))+2\nabla_2\nabla_1\xi(\bar x(0),\bar x(T))(V(\bar x(0)),
\\[2mm]& V(\bar x(T)))+\nabla_2^2\xi(\bar x(0),\bar x(T))(V(\bar x(T)),V(\bar x(T)))
\end{array}
\end{align}
for all $V\in\mathcal X(M)$.

For   $x,y\in M$  and a  differentiable curve $\gamma(\cdot):[a,b]\to M$ with $a,b\in \mathbb R$ and $a<b$ such that $\gamma(a)=x$ and $\gamma(b)=y$, we denote by $L_{xy}^\gamma$ the parallel translation along the curve $\gamma$. Especially, when $\rho(x,y)<\min\{i(x),i(y)\}$, there exists a unique shortest geodesic connecting $x$ and $y$, and we denote by $L_{xy}$ the parallel translation along this geodesic, from $x$ to $y$. We refer the readers  to \cite[Section 2.2]{cdz} for details of this notion.

Let $R$ denote the curvature tensor (see \cite[p.79]{pth}). For a vector $X\in T_xM$ and a covector $\eta\in T_x^*M$, with $x\in M$, we denote by $\tilde X\in T_x^*M$ and $\tilde\eta\in T_xM$ the dual covecor of $X$ and the dual vector of $\eta$, respectively,  i.e. 
\begin{align}\label{dcv}
\tilde X(Y)=\langle X,Y\rangle,\quad \langle\tilde\eta,Y\rangle=\eta(Y),\quad\forall\,Y\in T_xM.
\end{align}

Second, we recall some notions in tangent sets. For more details, we refer to \cite{AF}.
Let $X$ be a normed vector space with norm $|\cdot|_X$. Consider a subset $B\subset X$. 
The distance between a point $x\in X$ and $B$ is defined by $\textrm{dist}_B(x)=\inf\{|x-b|_X;\,b\in B\}$. If $\{B_h\}_{h>0}$ is a family of subsets of $X$, the lower limit of $\{B_h\}_{h>0}$ is given by $\lim\inf_{h\to0^+}B_h\stackrel{\triangle}{=}\{y\in X|\,\lim_{h\to0^+}\textrm{dist}_{B_h}(y)=0\}$. Denote by $\textrm{cl}\,B$ the closure of $B$. The adjacent cone to $B$ at a point $x\in\textrm{cl}\,B$ is defined by $T_B^\flat(x)=\lim\inf_{h\to0^+}\frac{B-\{x\}}{h}$ (see \cite[p. 127]{AF}).  Furthermore,  for $y\in T_B^\flat(x)$,  the second-order adjacent subset to $B$ at $(x,y)$ is given by $T_B^{\flat(2)}(x,y)=\lim\inf_{h\to0^+}\frac{B-x-hy}{h^2}$ (see \cite[Definition 4.7.2, p. 171]{AF}).

By the definition of the lower limit of a family of subset, one can represent $T_B^\flat(x)$ and $T_B^{\flat(2)}(x,y)$ (with $y\in T_B^\flat(x)$) as follows:
\begin{align*}
T_B^\flat(x)
=\{v\in X|\;\forall h_k\to 0^+,\,\exists \{v_k\}_{k=1}^{+\infty}\subseteq X\;s.t.\;\lim_{k\to+\infty}v_k= v\;\textrm{and}\;\forall k,\;x+h_kv_k\in B\},
\end{align*}
and
\begin{align*}
T_B^{\flat(2)}(x,y)=\left\{w\in X\left|\begin{array}{l}\forall h_k\to 0^+,\,\exists \{w_k\}_{k=1}^{+\infty}\subseteq X\;s.t.\;\lim_{k\to+\infty}w_k= w
\\ \textrm{and}\;\forall k,\;x+h_kv+h_k^2 w_k\in B
\end{array}\right.\right\}.
\end{align*}


Finally, the main conditions are stated as follows: 
\begin{description}
\item[($A1$)]  The map $f: [0,+\infty)\times M\times \mathbb R^m \to T M$, defined as $f(t,x,u)$, is measurable in $t$ and of class $C^1$ with respect to  $(x,u)$.
The map
$$\Phi(x_1,x_2)\equiv\big(\phi_0(x_1,x_2)^\top,
\\ \phi(x_1,x_2)^\top,\psi(x_1,x_2)^\top\big)^\top,\quad  x_1,x_2\in M
$$ 
 is $C^1$, and   there exists a constant $K>1$ and $x_0\in M$ such that,
\begin{equation}\begin{array}{l}\label{10}
\displaystyle |L_{x\hat x}f(t,x,u)-f( t,\hat x,u)|\leq \rho(x,\hat x),\;a.e.\,t\in[0,+\infty),\,u\in U,
\\[2mm]\displaystyle |\Phi(x_1,x_2)-\Phi(\hat x_1,\hat x_2)|_{\mathbb R^{r+j+k}}\leq  K(\rho(x_1,\hat x_1)+\rho(x_2,\hat x_2)),\;x_1,\hat x_1,x_2,\hat x_2\in M,
\\[2mm] \displaystyle |f(t,x_0,u)|\leq K\;a.e.\,t\in[0,+\infty),\;u\in U,\end{array}
\end{equation}  hold for   $x,\hat x\in M$ with $\rho(x,\hat x)< \min\{i(x),i(\hat x)\}$,  where $|\cdot|_{\mathbb R^p}$ (with $p\in\mathbb N$) denotes the norm in the space $\mathbb R^p$. Moreover, for each $r>0$, there exists a positive constant $K_r$ such that
\begin{align}\label{bdu}
|\nabla_uf(t,x,u)|\leq K_r,\quad\forall\,u\in U,\,x\in M\,\textrm{with}\,\rho(x_0,x)\leq r,\,\text{a.e. }\,t\in(0,+\infty),
\end{align}
 where $\nabla_u f(t,x,u)$ is defined by 
\begin{align}\label{ndu}
\nabla_u f(t,x,u)(\eta, V)\stackrel{\triangle}{=}\frac{\partial}{\partial\tau}\Big|_{\tau=0}f(t,x,u+\tau V)(\eta),
\end{align}
 for $V\in\mathbb R^m$ and $\eta\in T_x^*M$, and the corresponding norm is defined by
\begin{align*}
|\nabla_uf(t,x,u)|=\sup\{\nabla_u f(t,x,u)(\eta, V)|\,(\eta, V)\in T_x^*M\times\mathbb R^m,\,|\eta|\leq 1,\,|V|_{\mathbb R^m}\leq 1\}.
\end{align*}

\item[($A2$)]  The map $f(t,\cdot,\cdot)$ is $C^2$ for almost every    $t\in[0,+\infty)$,  where $\nabla_x^2f(t,x,u)$ (with $u\in U$) is the covariant differential of tensor $\nabla_xf(t,x,u)$ with respect to the variable $x$. The map $\Phi$
 is $C^2$ on $M\times M$, and   there exists a constant $K>1$ such that
\begin{align}\label{a2}\begin{array}{ll}
& \displaystyle |\nabla_xf(t,x_1,u)-L_{\hat x_1x_1}\nabla_x f(t,\hat x_1,u)|\leq K\rho(x_1,\hat x_1),
\\[2mm]& \displaystyle |\nabla_i\Phi(x_1,x_2)-L_{\hat x_ix_i}\nabla_i\Phi(\hat x_1,\hat x_2)|\leq K\left(\rho(x_1,\hat x_1)+\rho(x_2,\hat x_2)\right),\quad i=1,2,
\end{array}
\end{align}
hold for all $x_1,\hat x_1,x_2,\hat x_2\in M$ with $\rho(x_i,\hat x_i)<\min\{i(x_i),i(\hat x_i)\}$ ($i=1,2$), $u \in \mathbb R^m$ and almost every $t\in[0,+\infty)$, where   $\nabla_i\Phi$ and $|\nabla_i\Phi|$ (with $i=1,2$) are respectively given by (\ref{fdv}) and (\ref{cdv}). Moreover, for every $r>0$, there exists a positive constant $K_r>0$ such that
\begin{align}\label{sdub}
\max\{|\nabla_u\nabla_xf(t,x,u)|, |\nabla_u^2f(t,x,u)|\} \leq K_r,\quad\forall\,u\in U,\,x\in M\,\textrm{with}\,\rho(x_0,x)\leq r
\end{align}
holds for almost every $t\in(0,+\infty)$, where $\nabla_u^2f(t,x,u)$ and $\nabla_u\nabla_x f(t,x,u)$ (with $(t,x,u)\in[0,+\infty)\times M\times U$) are respectively defined by
\begin{align}\label{dfusx}\begin{array}{ll}
&\displaystyle \nabla_u^2f(t,x,u)(\eta,V,W)=\frac{\partial}{\partial\tau}\Big|_{\tau=0}\nabla_uf(t,x,u+\tau W)(\eta,V),
\\[2mm]&\displaystyle \nabla_u\nabla_xf(t,x,u)(\mu,X, V)=\frac{\partial}{\partial\tau}\Big|_{\tau=0}\nabla_xf(t,x,u+\tau V)(\mu, X),
\end{array}
\end{align}
for $(\mu, X,V,W)\in T_x^*M\times T_xM\times\mathbb R^m\times\mathbb R^m$, with the corresponding norms
\begin{align*}
&|\nabla_u^2f(t,x,u)|=\sup\{\nabla_u^2f(t,x,u)(\eta,V,W)|\,(\eta,V,W)\in T_x^*M\times\mathbb R^m\times\mathbb R^m,
\\&\quad\quad\quad\quad\quad\quad\quad \max\{|\eta|,|V|_{\mathbb R^m},|W|_{\mathbb R^m}\}\leq 1\},
\\& |\nabla_u\nabla_xf(t,x,u)|=\sup\{\nabla_u\nabla_x f(t,x,u)(\eta,X,W)|\,(\eta,X,W)\in T_x^*M\times T_xM\times\mathbb R^m,
\\ &\quad\quad\quad\quad\quad\quad\quad\quad\max\{|\eta|,|X|,|W|_{\mathbb R^m}\}\leq 1\}.
\end{align*}
\end{description}

\subsection{Necessary conditions for problem $(MP)_{r,T}$}

In  this section, we assume condition   $(A1)$ holds.
Let $(\bar u(\cdot),\bar x(\cdot))$ be a weak Pareto optimal solution of problem $(MP)_{r,T}$. For any map $\xi$ defined on $[0,T]\times M\times U$, we adopt the notation
\begin{align}\label{abbre}
\xi(t,\bar x(t),\bar u(t))\equiv\xi[t],\quad a.e.\,t\in[0,T]
\end{align}
for the sake of brevity.

In order to  formulate the necessary optimality conditions, we need to introduce  the Hamiltonian and Lagrangian functions. The map $H: [0,T]\times T^*M\times U\to\mathbb R$, given by 
\begin{align}\label{hcm}
H(t,x,\eta,u)\equiv \eta (f(t,x,u)),\quad (t,x,\eta,u)\in  [0,T]\times T^*M\times U,
\end{align}
is referred to as the Hamiltonian function.
Furthermore, for a vector $\ell\in \mathbb R^{r+j+k}$, we define the Lagrangian function  $\mathcal L_\ell: M\times M\to\mathbb R$   as follows:
\begin{align}\label{lfc}
\mathcal L_\ell(x_1,x_2)=\big(\phi_0(x_1,x_2)^\top,\phi(x_1,x_2)^\top,\psi(x_1,x_2)^\top\big)\ell,\quad\forall\,x_1,x_2\in M.
\end{align}
For this vector $\ell$, we also define a covector  field $p^\ell(\cdot)$ along $\bar x(\cdot)$, i.e. $p^\ell(t)\in T_{\bar x(t)}^*M$ for each $t\in[0,T]$, as the solution to the following differential equation:
\begin{align}\label{de}
\left\{\begin{array}{l}
\nabla_{\dot{\bar x}(t)} p^\ell=-\nabla_xf[t](p^\ell(t),\cdot),\quad a.e.\,t\in(0,T),
\\ p^\ell(0)=-d_1\mathcal L_\ell(\bar x(0),\bar x(T)),\quad p^\ell(T)=d_2\mathcal L_\ell(\bar x(0),\bar x(T)),
\end{array}\right.
\end{align}
where $[t]$ is defined by (\ref{abbre}), and $d_i\mathcal L_\ell(x_1,x_2)$ (for $i=1,2$) denotes  the differential of $\mathcal L_\ell$ with respect to $x_i$.

The first-order necessary optimality condition for a weak Pareto optimal solution is given below.
\begin{theorem}\label{fon}
Assume that  the set $U$ is closed and convex,  condition  $(A1)$  holds, and  $(\bar u(\cdot),\bar x(\cdot)) $  is a weak Pareto optimal solution of problem $(MP)_{r,T}$, where $\bar u(\cdot)\in L^2(0,T;\mathbb R^m)\cap \mathcal U_T$. Set
$
 I_A=\{0_1,\cdots,0_r\}\cup\{i\in\{1,\cdots,j\}|\;\phi_i(\bar x(0),\bar x(T))=0\}$. 
 Then, there exists a nonzero vector 
\begin{align}\label{lm1}
\ell=(\ell_{0_1},\cdots,\ell_{0_r},\ell_1,\cdots,\ell_j,\ell_\psi^\top)^\top\in\mathbb R^{r+j+k}
\end{align}
satisfying
\begin{align}\label{lm2}
\ell_i\leq 0\;\; \textrm{if}\;\; i\in I_A;\quad\ell_i=0\;\;\textrm{if}\;\;i\in \{1, \cdots,j\}\setminus I_A,
\end{align}
 such that
\begin{align}\label{lm3}
\nabla_uH(t,p^\ell(t),\bar x(t),\bar u(t))(v(t))\leq 0,\quad a.e.\,t\in(0,T)
\end{align}
holds for all $v(\cdot)\in L^2(0,T;\mathbb R^m)$ with 
$
v(t)\in T_U^\flat(\bar u(t))
$ for a.e. $t\in(0,T)$, 
where $p^\ell(\cdot)$ is the solution to (\ref{de}), and $\nabla_u H(t,x,p,u)$,  with $(t,x,p,u)\in[0,T]\times T^*M\times U$, represents the partial derivative of $H$ with respect to  
 $u$.
\end{theorem} 

Based on the first-order necessary optimality condition presented in Theorem \ref{fon}, we introduce the Lagrangian multiplier and singular direction as follows.
\begin{definition}\label{lm}
Assume  the assumptions in Theorem \ref{fon} hold.  Denote by $O_{\mathbb R^{r+j+k}}$ the origin of the space $\mathbb R^{r+j+k}$. A vector $\ell\in\mathbb R^{r+j+k}\setminus\{O_{\mathbb R^{r+j+k}}\}$ is called a Lagrangian multiplier  of problem $(MP)_{r,T}$, if it satisfies (\ref{lm1}) and (\ref{lm2}), and  there is a covector field $p^\ell(\cdot)$ satisfying (\ref{de}) and (\ref{lm3}).
\end{definition}

\begin{definition}\label{sigd}
Assume the assumptions in Theorem \ref{fon} hold. A function  $v(\cdot)\in L^2(0,T;\mathbb R^m)$ is called a singular direction, if there exists a vector field $X_v(\cdot)$ along $\bar x(\cdot)$, i.e. $X_v(t)\in T_{\bar x(t)} M$ for each $t\in[0,T]$, satisfying
\begin{align}\label{sdc}
\nabla_{\dot{\bar x}(t)}X_v=\nabla_xf[t](\cdot, X_v(t))+\nabla_uf[t](\cdot,v(t)),\quad a.e.\,t\in[0,T],
\end{align}
and
\begin{align}\label{sdc1}\begin{array}{ll}
&v(t)\in T_U^\flat(\bar u(t)),\quad a.e.\,t\in[0,T],
\\  & \nabla_1 \phi_i(\bar{x}(0), \bar{x}(T))\left(X_v(0)\right)+\nabla_2 \phi_i(\bar{x}(0), \bar{x}(T))\left(X_v(T)\right) \leq 0, \forall i \in I_A, 
\\ & \nabla_1 \psi(\bar{x}(0), \bar{x}(T))\left(X_v(0)\right)+\nabla_2 \psi(\bar{x}(0), \bar{x}(T))\left(X_v(T)\right)=0,
\end{array}
\end{align}
where   $\nabla_i\psi$ is defined by (\ref{fdv}), the notation $[t]$ is given by (\ref{abbre}), and for $(t,u)\in[0,T]\times U$,   $\nabla_xf(t,\cdot,u)$ represents the  covariant derivative of tensor field $f(t,\cdot,u)$  with respect to $x$, while $\nabla_uf$ is defined by (\ref{ndu}).
\end{definition}
 
 The following result reveals what happens along a singular direction.
\begin{proposition}\label{psd}
Assume the assumptions in Theorem \ref{fon} hold. If $v(\cdot)$  is a singular direction, then for any Lagrangian multiplier of problem $(MP)_{r,T}$, it holds that
\begin{align*}
\nabla_u H(t, p^\ell(t),\bar x(t),\bar u(t))(v(t))=0,\quad a.e.\,t\in(0,T).
\end{align*}
\end{proposition}
 The proof of the above result will be given in Section \ref{ae21}. 

Regarding the first-order necessary condition and singular directions, we have a comment below.
\begin{rem}
Assume the conditions in Theorem \ref{fon} hold. The first-order necessary condition can be understood as the existence of a Lagrangian multiplier. We can conclude from Proposition \ref{psd} and (\ref{lm3}) that, along a singular direction, the first-order necessary condition degenerates. Thus, we need to investigate the second-order necessary conditions along the singular direction.
\end{rem}

For a singular direction $v(\cdot)\in L^2(0,T;\mathbb R^m)$,  we set 
\begin{align*}\begin{array}{l}
I_A^\prime\stackrel{\triangle}{=}\Big\{i\in I_A|\,\nabla_1\phi_i(\bar x(0),\bar x(T) )(X_v(0))+\nabla_2\phi_i(\bar x(0),\bar x(T))(X_v(T))=0\Big\},
\end{array}
\end{align*}
where $X_v(\cdot)$ is given in  Definition  \ref{sigd}. The set
$I_A^\prime$ can be understood as an index set that contains all indices of the components for which both the zeroth-order and first-order constraints are active.

Moreover, if assumptions $(A1)$ and $(A2)$ hold,  the map $\nabla_u\nabla_x H(t,x,p,u): T_x M\times\mathbb R^m\to\mathbb R$, with $(t,x,p,u)\in[0,T]\times T^*M\times\mathbb R^m$, is defined by
\begin{align*}
\nabla_u\nabla_x H(t,x,p,u)(X,V)=\lim_{\epsilon\to 0}\frac{\nabla_x H(t,x,p,u+\epsilon V)(X)-\nabla_x H(t,x,p,u)(X)}{\epsilon}, 
\end{align*}
for all $(X,V)\in T_xM\times\mathbb R^m$.  Additionally, $\nabla_u^2 H(t,x,p,u)$ denotes the Hessian of $H$ with respect to $u$.

The second-order necessary optimality condition is stated as follows.
\begin{theorem}\label{sen}
Assume that   the set $U$ is convex and closed, conditions $(A1)$ and $(A2)$ hold,  and $(\bar u(\cdot),\bar x(\cdot))$ is a weak Pareto optimal solution of problem $(MP)_{r,T}$, where $\bar u(\cdot)\in L^2(0,T;\mathbb R^m)$. Let $v(\cdot)\in L^2(0,T;\mathbb R^m)$ be a singular direction, with $X_v(\cdot)$ satisfying (\ref{sdc}) and (\ref{sdc1}). Assume  there exists $\eta(\cdot)\in L^2(0,T;\mathbb R^m)$ and $\alpha>0$ such that $d_U\left(\bar u(t)+\epsilon v(t)\right)\leq \epsilon^2 \eta(t)$ holds for all $\epsilon\in[0,\alpha]$ and almost every $t\,\in[0,T]$, 
and  the set $\mathcal S_T\stackrel{\triangle}{=}\{\sigma(\cdot)\in L^2(0,T;\mathbb R^m)|\,\,\sigma(t)\in T_U^{\flat(2)}(\bar u(t),v(t))\;a.e.\,t\in[0,T]\}\not=\emptyset$. Then,
 there exists a Lagrangian multiplier for problem $(MP)_{r,T}$, denoted by   
 $\ell=(\ell_{0_1},\cdots,\ell_{0_r},\ell_1,\cdots,\ell_j,\ell_\psi^\top)^\top$, which satisfies
\begin{align}\label{slm}
 \ell_i=0,\;\textrm{if}\;i\notin I_A^\prime,
\end{align}
such that the inequality
\begin{align}\label{sonc2}\begin{array}{l}
\displaystyle\int_0^T\Big(2\nabla_u H\{t;\ell\}\left(\sigma(t)\right)+ \nabla_x^2 H\{t;\ell\}\left( X_v(t),X_v(t)  \right)+2\nabla_u\nabla_x H\{t;\ell\}\left( X_v(t),v(t) \right)
\\[2mm]\quad\quad\displaystyle+ \nabla_u^2 H\{t;\ell\}(v(t),v(t))- R\left( \tilde p^\ell(t),X_v(t),f[t], X_v(t)  \right)\Big)dt
+ D^2\mathcal L_\ell(\bar x(0),\bar x(T))(X_v)
\leq 0
\end{array}
\end{align}
holds for all $\sigma(\cdot)\in\mathcal S_T$. Here, $p^\ell(\cdot)$ is the solution to (\ref{de}),   $\tilde p^\ell(t)$ is the dual vector of $p^\ell(t)$ as defined by (\ref{dcv}),  and  $D^2\mathcal L_\ell(\bar x(0),\bar x(T))(X_u)$ is given by  (\ref{msd}). Additionally,   for any $t\in[0,T]$, the notation $\{t;\ell\}$ abbreviates $(t,\bar x(t), p^\ell(t),\bar u(t)))$.
\end{theorem}

\subsection{Necessary conditions for problem $(MP)_r$}
We need to introduce additional conditions:
\begin{description}
\item[$(A1^\prime)$] The maps $f$ and $\Phi$ introduced in assumption $(A1)$ are $C^1$, and there exists a constant $K\geq 1$ and $x_0\in M$ such that the second line of (\ref{10}) holds, and
\begin{align}\label{10p}\begin{array}{l}
\displaystyle |L_{x\hat x}f(t,x,u)-f(\hat t,\hat x,u)|\leq K\big(|t-\hat t|+\rho(x,\hat x)\big), \quad\forall\,u\in U,\;t,\hat t\in[0,+\infty),
\\ |f(0,x_0,u)|\leq K,\quad\forall\,u\in U,
\end{array}
\end{align}
hold for     $x,\hat x\in M$ with $\rho(x,\hat x)\leq \min\{i(x),i(\hat x)\}$.  Moreover, for each $r>0$, there exists $K_r>0$ such that 
\begin{align*}
|\nabla_uf(t,x,u)|\leq K_r,\;\forall\,(t,x,u)\;\in[0,+\infty)\times M\times U\;\textrm{with}\; |t|+\rho(x,x_0)\leq r,
\end{align*}
where $\nabla_u f(t,x,u)$ (with $(t,x,u)\in[0,+\infty)\times  M\times U$) is given by (\ref{ndu}).

\item[($A2^\prime$)]  The maps $f$ and $\Phi$ are $C^2$. Furthermore, there exists a constant $K>1$ such that the second line of (\ref{a2}) holds for all $x_1,x_2,\hat x_1,\hat x_2\in M$ with $\rho(x_l,\hat x_l)<\min\{i(x_l),i(\hat x_l)\}$ ($l=1,2$), and the following inequalities
\begin{align*}
&|f_t(t,x,u)-L_{\hat x x}f_t(\hat t,\hat x,u)|\leq K(|t-\hat t|+\rho(x,\hat x)),\;\forall\,t,\hat t \in[0,+\infty),\, u\in U,
\\ & |\nabla_xf(t,x,u)-L_{\hat x x}\nabla_xf(\hat t,\hat x,u)|\leq K(|t-\hat t|+\rho(x,\hat x)),\;\;\forall\,t,\hat t \in[0,+\infty),\, u\in U,
\end{align*}
hold for all  $x,\hat x\in M$ with $\rho(x,\hat x)<\min\{i(x),i(\hat x)\}$, where $ f_t(t,x,u)$ is defined as
$
f_t(t,x,u)=\lim_{\epsilon\to0}\frac{f( t+\epsilon,x,u )-f(t,x,u)}{\epsilon}$. 
Moreover, for every $r>0$, there exists a positive constant $K_r$ such that
\begin{align*}
&\max\{|f_t(t,x,u)|, |\nabla_xf(t,x,u)|, |\nabla_uf_t(t,x,u)|, |\nabla_u^2f(t,x,u)|,
|\nabla_u\nabla_xf(t,x,u) |\}\leq K_r
\end{align*}
holds for all $u\in U$ and $(t,x)\in[0,+\infty)\times M$ with $|t|+\rho(x,x_0)\leq r$.

\end{description}

The first-order necessary optimality condition is stated as follows.

\begin{theorem}\label{fnct} Assume that condition $(A1)$  holds,   $U$ is closed and convex, and $(\bar u(\cdot),\bar x(\cdot), \bar T)$ is a weak Pareto optimal solution with $\bar u(\cdot)\in L^2_{loc}(0,+\infty;\mathbb R^m)\cap \mathcal U$.  
Set
$\hat I_A=\{0_1,\cdots,0_r\}\cup\{i\in\{1,\cdots,j\}\mid \phi_i(\bar{x}(0),\bar{x}(\bar T))=0\}$.
Let $H$ and $\mathcal L_\ell$ be given by (\ref{hcm}) and (\ref{lfc}), respectively.
 Then, there exists a nonzero vector
$
\ell=(\ell_{{0_{1}}},\cdots,\ell_{{0_{r}}},\ell_{1},\cdots,\ell_{j},\ell_{\psi}^{\top})^{\top}\in\mathbb{R}^{r+j+k}
$
satisfying
\begin{align}\label{aim}
\ell_{i}\leq0\;\textrm{if}\;i\in \hat I_{A};\quad\ell_{i}=0\;\textrm{if}\;i\in\{1,\cdots,j\}\setminus \hat I_{A},
\end{align}
such that 
\begin{align}\label{font}\begin{array}{l}
H\left(t,\bar x(t),p^\ell(t),\bar u(t)   \right)+\int_t^{\bar T}\nabla_t H\left(\mu,\bar x(\mu), p^\ell(\mu),\bar u(\mu)\right)d\mu =0,\quad a.e.\,t\in[0,\bar T],
\\[2mm] \nabla_u H\left(t,\bar x(t),p^\ell(t),\bar u(t)   \right)(v(t))\leq0,\quad a.e.\,t\in[0,\bar T],
\end{array}
\end{align}
hold for all $v(\cdot)\in L^2(0,\bar T;\mathbb R^m)$ with $v(t)\in T_U^\flat(\bar u(t))$ a.e. $t\in[0,\bar T]$, where $p^\ell(\cdot)$ satisfies
\begin{align}\label{detd}
\left\{\begin{array}
{l} \nabla_{\dot{\bar x}(t)} p^\ell=-\nabla_xf[t](p^\ell(t),\cdot),\quad a.e.\,t\in(0,\bar T),
\\ p^\ell(0)=-d_1\mathcal{L}_\ell\left(\bar{x}(0),\bar{x}(\bar T)\right), p^\ell(\bar T)=d_2\mathcal{L}_\ell\left(\bar{x}(0),\bar{x}(\bar T)\right),
\end{array}
\right.
\end{align}
and $\nabla_t H(t,x,p,u)$, with $(t,x,p,u)\in [0,+\infty)\times T^*M\times \mathbb R^m$, denotes the partial derivative of $H$ with respect to $t$.
\end{theorem}

Similar to the case of problem $(MP)_{r,T}$, we introduce the notion of  Lagrangian multiplier.

\begin{definition}\label{lmtd} Assume  the assumptions in Theorem \ref{fnct} hold. A nonzero vector $\ell\in\mathbb R^{r+j+k} $ is called a Lagrangian multiplier for problem $(MP)_r$,  if it satisfies (\ref{aim}), (\ref{font}) and (\ref{detd}). Furthermore,  a pair $(\xi(\cdot), v(\cdot))\in L^2(0,\bar T;\mathbb R\times\mathbb R^m)$ is called a singular direction,  if there exists a vector field $X_{\xi,v}(\cdot)$ along $\bar x(\cdot)$ satisfying
\begin{align}\label{foet}\begin{array}{ll}
\displaystyle\nabla_{\dot{\bar  x}(t)} X_{\xi,v}=&\nabla_xf[t](\cdot, X_{\xi,v}(t))+\frac{1}{\bar T}\int_0^t\xi(s)ds\,\nabla_tf[t]+\frac{\xi(t)}{\bar T}f[t]
\\[2mm]&\displaystyle+\nabla_u f[t](v(t)),\quad a.e.\, t\in(0,\bar T),
\end{array}
\end{align}
and
\begin{align}\label{sdt}\begin{array}{ll}
&\displaystyle v(t)\in T_U^\flat(\bar u(t)),\quad a.e.\,t\in(0,\bar T),
\\[2mm]& \nabla_1\phi_i(\bar x(0),\bar x(\bar T))X_{\xi,v}(0)+\nabla_2\phi_i(\bar x(0),\bar x(\bar T))X_{\xi,v}(\bar T)\leq 0,\quad i\in \hat I_A,
\\[2mm]&\displaystyle\nabla_1\psi(\bar x(0),\bar x(\bar T))X_{\xi,v}(0)+\nabla_2\psi(\bar x(0),\bar x(\bar T))X_{\xi,v}(\bar T)=0.
\end{array}
\end{align}
With respect to $(\xi(\cdot), v(\cdot))$, we set
\begin{align*}
{\hat I}^\prime_A\triangleq\Big\{i\in \hat I_A|\,\nabla_1\phi_i(\bar x(0),\bar x(\bar T))(X_{\xi,v}(0))+\nabla_2\phi_i(\bar x(0),\bar x(\bar T))(X_{\xi,v}(T))=0\Big\}.
\end{align*}
\end{definition}

The second-order necessary optimality condition is as follows.

\begin{theorem}\label{snt}
Assume that   conditions in Theorem \ref{fnct} and assumption $(A2^\prime)$ hold. 
 Let $(\xi(\cdot),v(\cdot))\in L^2(0,\bar T;\mathbb R\times \mathbb R^m)$ be a singular direction, with  $X_{\xi,v}(\cdot)$ satisfying  (\ref{foet}) and (\ref{sdt}). We also assume that
 $\mathcal S_{\bar T}\stackrel{\triangle}{=}\{\sigma(\cdot)\in L^2(0,\bar T;\mathbb R^m)|\,\sigma(t)\in T_U^{\flat(2)}(\bar u(t),v(t))\;a.e.\,t\in[0,\bar T]\}\not=\emptyset.$  Then, there exists a Lagrangian multiplier  $\ell=(\ell_{0_1},\cdots,\ell_{0_r},\ell_1,\cdots,\ell_j,\ell_\psi^\top)^\top$,
satisfying
\begin{align}\label{snti}
\ell_i=0,\quad\textrm{if}\;i\not\in{\hat I}^\prime_A,
\end{align}
 such that the following inequality 
\begin{align}\label{snct1}\begin{array}{l}
\displaystyle\int_0^{\bar T}\Big(2\nabla_u H\{t;\ell\}\left(\sigma(t)\right)+ \nabla_x^2 H\{t;\ell\}\left( X_{\xi,v}(t),X_{\xi,v}(t)  \right)+2\nabla_u\nabla_x H\{t;\ell\}\left( X_{\xi,v}(t),v(t) \right)
\\[2mm] \displaystyle+ \nabla_u^2 H\{t;\ell\}(v(t),v(t))- R\left( \tilde p^\ell(t),X_{\xi,v}(t),f[t], X_{\xi,v}(t)  \right)\Big)dt
\\ [2mm]\displaystyle+\int_0^{\bar T}\Bigg(\nabla_t^2H\{t;\ell\}\left(\frac{1}{\bar T}\int_0^t\xi(\mu)d\mu  \right)^2+2\nabla_x\nabla_t H\{t;\ell\}\left( \frac{1}{\bar T}\int_0^t\xi(\mu)d\mu,X_{\xi,v}(t) \right)
\\[2mm] \displaystyle +\frac{2}{\bar T^2}\int_0^t\xi(\mu)d\mu\xi(t)\nabla_t H\{t;\ell\}+2\nabla_u\nabla_tH\{t;\ell\}\left(\frac{1}{\bar T}\int_0^t\xi(\mu)d\mu,v(t)  \right)
\\[2mm] \displaystyle +\frac{2}{\bar T}\xi(t)\nabla_xH\{t;\ell\}(X_{\xi,v}(t))
+\frac{2}{\bar T}\xi(t)\nabla_u H\{t;\ell\}(v(t))\Bigg)dt
+ D^2\mathcal L_\ell(\bar x(0),\bar x(\bar T))(X_{\xi,v})
\leq 0
\end{array}
\end{align}
holds for all $\sigma(\cdot)\in\mathcal S_{\bar T}$, where,  for $t\in[0,\bar T]$, the notation $\{t;\ell\}$ denotes $(t, \bar x(t), p^\ell(t), \bar u(t))$, with $p^\ell(\cdot)$ being the solution to  (\ref{detd}).
\end{theorem}

As an application of the main results, we provide an example below.

\begin{example}\label{exg} Set
\begin{align*}
f(x,u)
& \stackrel{\triangle}{=}\left( \begin{array}{c}u_2 x_3^2
\\h(x_1,x_2,u_1,u_2)
\\[2mm] \displaystyle\frac{2u_2x_1x_3^2+2x_2h(x_1,x_2,u_1,u_2)}{1+x_1^2+x_2^2}
\end{array} \right),\;(x,u)=(x_1,x_2,x_3,u_1,u_2)\in\mathbb R^3\times\mathbb R^2,
\end{align*}
where $h(x_1,x_2,u_1,u_2)=-x_1^2+4x_1u_2-u_1$
for $(x_1,x_2,u_1,u_2)\in\mathbb R^2\times \mathbb R^2$.
Given $T>0$, consider the following optimal control problem:
Minimize  $(-x_1(T)^2,-x_3(T))$ over $(x_1(\cdot), x_2(\cdot), x_3(\cdot), u_1(\cdot), u_2(\cdot))$
subject to
 \begin{align}\label{exs}
\left\{\begin{array}{l} \dot x(t) =f(x(t),u(t)),
\quad  a.e.\,t\in(0,T),
\\[2mm]\displaystyle x_1(0)=0,\quad x_2(0)\leq 0,\quad x_1(T)^3+ x_2(T)+T=0,
\\[2mm]\displaystyle x_3(t)=\ln(1+x_1(t)^2+x_2(t)^2),\quad \forall\, t\in[0,T],
\\[2mm]\displaystyle (u_1(t),u_2(t))\in B_{\mathbb R^2}(O_{\mathbb R^2},1),\quad a.e.\,t\in[0,T], \end{array}\right.
\end{align}
where  $O_{\mathbb R^2}$ is the origin of  space $\mathbb R^2$, and $B_{\mathbb R^2}(O_{\mathbb R^2},1)=\{y\in\mathbb R^2|\;|y|\leq 1\}$. Please show that the admissible control $\bar u(\cdot)=(\bar u_1(\cdot),\bar u_2(\cdot))=(1,0)$ is not a weak Pareto optimal solution.

\end{example}
Before solving this example, we have a remark.
\begin{rem}\label{reer}Regarding the pointwise state constraint ``\( x_3(t) = \ln(1 + x_1(t)^2 + x_2(t)^2) \) for all \( t \in [0, T] \)'', we have two comments. First, without this constraint, the corresponding control system may not necessarily have a solution on \( [0, T] \). This is because \( f(x, u) \) is only locally Lipschitz continuous with respect to the variable \( x \in \mathbb{R}^3 \). According to the existing theory of ordinary differential equations, this cannot guarantee the existence of a solution. However, the set of pointwise state constraints forms a Riemannian manifold. On this manifold, \( f(x, u) \) is Lipschitz continuous with respect to the variable \( x \) under the Riemannian metric. Consequently, by the theory of ordinary differential equations on Riemannian manifolds, the solution to the control system (\ref{exs}) exists. Second, to the best of our knowledge, the presence of this pointwise state constraint and the  equality-type terminal time constraint ``$x_1(T)^3+ x_2(T)+T=0$'' renders existing results on optimality conditions (such as those in the literature \cite{kwy, BA04, zhu00, kty, KQWY20, kb23}) inapplicable to this problem.
\end{rem}

{\it Solution of  Example \ref{exg}.}\;  First, we demonstrate that the control system (\ref{exs}) evolves on the manifold 
$ M = \{(x_1,x_2,a(x_1,x_2)) | (x_1,x_2) \in \mathbb{R}^2\}, $
where $ a(x_1,x_2) = \ln(1 + x_1^2 + x_2^2) $ for $ (x_1,x_2) \in \mathbb{R}^2 $. According to Proposition \ref{csg1}, we conclude that $ M $ is a two-dimensional Riemannian manifold equipped with the coordinate chart  (\ref{ccg}),  the basis (\ref{btg}) for the tangent space of $ M $ at $ (x_1,x_2,a(x_1,x_2)) $,
 the Riemannian metric   (\ref{rmg}),
and the sectional curvature   (\ref{scg}).
Given that
\begin{align}\label{lef2}f(x,u)
=\sum_{i=1}^2f_i(x,u)\frac{\partial}{\partial x_i}\Big|_M,
\end{align}
for $x=(x_1,x_2,x_3)\in M$ and $u=(u_1,u_2)\in\mathbb R^2$, where
$
f_1(x,u)=u_2x_3^2$ and $  f_2(x,u)=h(x_1,x_2,u)$. 
It follows that the control system (\ref{exs}) indeed evolves on the manifold $ M $.

Second, we verify that this optimal control problem satisfies conditions $(A1)$ and $(A2)$. Specifically, we need to show that (\ref{10}) and (\ref{a2}) hold. To achieve this, let us define the matrix $G$ as a function of $(x_1,x_2)\in\mathbb{R}^2$:
$
G(x_1,x_2)=(g_{il}(x_1,x_2))$, where $g_{il}(x_1,x_2)$ is defined by  (\ref{rmg}).

To prove the first line of (\ref{10}),  by \cite[Lemma 4.1]{cdz}, it suffices to show that  there exists a  positive constant $K$ such that
\begin{align}\label{bcfxu}
|\nabla_xf(x,u)|\leq K,\quad  \forall\,(x,u)\in M\times B_{\mathbb R^2}(O_{\mathbb R^2},1).
\end{align}
 
 We will prove this in the coordinate chart  $(M,\varphi)$. 
Assume that
$\nabla_xf(x,u)=\\  \sum_{i,l=1}^2a_{il}(x_1,x_2,u)\frac{\partial}{\partial x_i}\Big|_M\otimes dx_l|_M,$ where $\{dx_1|_M, dx_2|_M\}$ forms the dual basis of $\big\{ \frac{\partial}{\partial x_1}\big|_M, \frac{\partial}{\partial x_2}\big|_M\big\}$, meaning that
\begin{align}\label{dbcs}
dx_i|_M\left(\frac{\partial}{\partial x_l}\Big|_M\right)=\delta_i^l,\quad i,l=1,2,
\end{align}
with $\delta_i^l$ being Kronecker delta symbols.

 According to the definitions of the covariant differential (see \cite[Definition 5.7, p.102]{c}) and   Christoffel symbols (see \cite[p.55]{c}), as well as properties of the covariant derivative (see \cite[Theorem 2.2.2, p.53]{pth}), along with (\ref{lef2}) and (\ref{dbcs}),  we find that
\begin{align*}
&a_{11}(x_1,x_2,u)
\\&=\nabla_xf\left(  dx_1|_M,\frac{\partial}{\partial x_1}\Big|_M\right)
\\&=\nabla_{\frac{\partial}{\partial x_1}\Big|_M}\sum_{i=1}^2f_i(\cdot,u)\frac{\partial}{\partial x_i}\Big|_M\left(dx_1|_M \right)
\\&=\frac{\partial}{\partial x_1}\Big|_M f_1(\cdot, u)+f_1(x,u)\nabla_{\frac{\partial}{\partial x_1}\Big|_M}\frac{\partial}{\partial x_1}\Big|_M(dx_1|_M)+f_2(x,u)\nabla_{\frac{\partial}{\partial x_1}\Big|_M}\frac{\partial}{\partial x_2}\Big|_M(dx_1|_M)
\\&=\frac{4x_1a(x_1,x_2)u_2}{1+x_1^2+x_2^2}+f_1\left(x_1,x_2,a(x_1,x_2),u\right)\Gamma_{11}^1(x_1,x_2)+h(x_1,x_2,u)\Gamma_{12}^1(x_1,x_2),
\end{align*}
where the Christoffel symbols $\Gamma_{ip}^l$ (with $i,p,l=1,2$) are given by  (\ref{csg}). Similarly, we can obtain
\begin{align*}
a_{12}(x_1,x_2,u)=&f_1(x_1,x_2,a(x_1,x_2),u)\Gamma_{21}^1(x_1,x_2)+h(x_1,x_2,u)\Gamma_{22}^1(x_1,x_2),
\\a_{21}(x_1,x_2,u)=&f_1\left(x_1,x_2,a(x_1,x_2),u\right)\Gamma_{11}^2(x_1,x_2)+\frac{\partial}{\partial x_1}h(x_1,x_2,u)
+h(x_1,x_2,u)\Gamma_{12}^2(x_1,x_2),
\\a_{22}(x_1,x_2,u)=&f_1\left(x_1,x_2,a(x_1,x_2),u\right)\Gamma_{21}^2(x_1,x_2)+\frac{\partial}{\partial x_2}h(x_1,x_2,u)
+h(x_1,x_2,u)\Gamma_{22}^2(x_1,x_2).
\end{align*}
Let $A(x_1,x_2,u)=(a_{il}(x_1,x_2,u))$ for $(x_1,x_2,u)\in\mathbb R^2\times\mathbb R^m$.
Thus, using Proposition \ref{ntlc},  for $V=(V_1,V_2)\in\mathbb R^2$ and $X=\sum_{i=1}^2 X_i\frac{\partial}{\partial x_i}\Big|_M\in T_xM$, we get
\begin{align*}
\nabla_xf(x,u)(\eta,X)
=&\begin{pmatrix}  X_1&X_2
\end{pmatrix} A(x_1,x_2,u)\begin{pmatrix}
\eta_1&\eta_2
\end{pmatrix}^\top
\\\leq&|X|\left|\sqrt{G}^{-1}A(x_1,x_2,u)\sqrt{G}\right|_{\mathbb R^{2\times 2}}|\eta|,
\end{align*}
which implies
$
|\nabla_xf(x,u)|\leq  \left|\sqrt{G(x_1,x_2)}^{-1}A(x_1,x_2,u)\sqrt{G(x_1,x_2)}\right|_{\mathbb R^{2\times 2}},
$
where $|\cdot|_{\mathbb R^{2\times 2}}$ denotes the norm for a square matrix of order 2.

Similarly we can obtain $|\nabla_uf(x,u)|\leq |B(x_1,x_2,u)\sqrt{G(x_1,x_2)}|_{\mathbb R^{2\times 2}}$, where  $B(x_1,x_2,u)=\\\left(\frac{\partial}{\partial u_i}f_l(x_1,x_2,a(x_1,x_2),u)\right)_{i,l=1}^2$.

Given that the eigenvalues of $G(x_1,x_2)$ are
$
\lambda_1(x_1,x_2)=1+\frac{4(x_1^2+x_2^2)}{(1+x_1^2+x_2^2)^2}$ and  $\lambda_2(x_1,x_2)=1$, 
there exists an orthogonal matrix  $P(x_1,x_2)$ such that 
\begin{align*}
\sqrt{G(x_1,x_2)}=P(x_1,x_2)^\top\begin{pmatrix}\sqrt{\lambda_1(x_1,x_2)}&0\\0&\sqrt{\lambda_2(x_1,x_2)}\end{pmatrix}P(x_1,x_2),
\\ \left.\left.\left.\sqrt{G(x_1,x_2)}^{-1}=P(x_1,x_2)^\top\left(\begin{array}{cc}\frac{1}{\sqrt{\lambda_1(x_1,x_2)}}&0\\0&\frac{1}{\sqrt{\lambda_2(x_1,x_2)}}\end{array}\right.\right.\right.\right)P(x_1,x_2).
\end{align*}
Since  
$1\leq |\lambda_1(x_1,x_2)|\leq 5$ for all $(x_1,x_2)\in\mathbb R^2$, we have
$|\sqrt{G(x_1,x_2)}|_{\mathbb R^{2\times2}}\leq \sqrt 5$ and $|\sqrt{G(x_1,x_2)}^{-1}|_{\mathbb R^{2\times2}}\leq 1$ for all $(x_1,x_2)\in\mathbb R^2$.
Furthermore, based on the special structure of the Christoffel symbols (\ref{csg}),   there exists a positive constant $C$ such that  $\left|A(x_1,x_2,u)\right|\leq C$ holds for all 
$(x_1,x_2,u)\in\mathbb R^2\times  B_{\mathbb R^2}(O_{\mathbb R^2},1)$.
Therefore, there exists a positive constant $K$ such that    (\ref{bcfxu}) holds.  Moreover, for each $r>0$, there exists $K_r>0$ such that (\ref{bdu}) holds.

Using similar methods, we can prove that the second line of  (\ref{10}), the second line of   (\ref{a2}) and (\ref{sdub})  also hold. Consequently,  assumptions $(A1)$ and $(A2)$ are verified.

Finally, we  apply   Theorem \ref{fon} and Theorem \ref{sen} to this problem. 
 In the coordinate chart $(M,\varphi)$, according to (\ref{lef2}), the system (\ref{exs}) can be expressed as follows:
\begin{align*}
\left\{\begin{array}{l}
\sum_{i=1}^2\dot x_i(t)\frac{\partial}{\partial x_i}\Big|_M=\sum_{i=1}^2f_i\left(x_1(t), x_2(t), a(x_1(t),x_2(t)), u(t)\right)\frac{\partial}{\partial x_i}\Big|_M,
 \; a.e.\,t\in(0,T),
\\ \psi(x_1(0), x_1(T), x_2(T))=0,
\quad \phi(x_2(0))\leq 0;
\quad (u_1(t),u_2(t))\in  B_{\mathbb R^2}(O_{\mathbb R^2},1),\;a.e. \;t\in(0,T),
\end{array}
\right.
\end{align*}
where   $\psi(x_1(0),x_1(T), 
x_2(T))=(x_1(0),x_1(T)^3+x_2(T)+T)$ and $\phi(x_2(0))=x_2(0)$ for $(x_1(0),x_2(0), x_1(T),\\ x_2(T))\in\mathbb R^4$. Furthermore,  from Proposition \ref{elsn},  we obtain that
the Hamiltonian function   is given by the following expression:
$
H(x_1,x_2,p_1,p_2,u)=u_2a(x_1,x_2)^2p_1+h(x_1,x_2,u)p_2,
$
for all $(x_1,x_2,p_1,p_2,u)\in\mathbb R^2\times\mathbb R^2\times B_{\mathbb R^2}(O_{\mathbb R^2},1)$.

Assuming $\bar u(\cdot)=(1,0)$ is a weak Pareto optimal control,  the corresponding trajectory is $\bar x_1(\cdot)=0$, $\bar x_2(t)=-t$ and $\bar x_3(t)=\ln(1+t^2)$ for $t\in[0,T]$.  For $\bar u(\cdot)$, applying  Theorem  \ref{fon} and using Proposition \ref{elsn}, we find that, there exists 
a vector $\ell=(\ell_{0_1},\ell_{0_2},\ell_\phi,\ell_{\psi_1},\ell_{\psi_2})\in\mathbb R^5\setminus\{O_{\mathbb R^5}\}$ such that
$
\nabla_u H(\bar x_1(t),\bar x_2(t), p^\ell(t), \bar u(t))(v_1(t),v_2(t))\leq 0$  a.e. $t\in(0,T)
$
holds for all $(v_1(\cdot),v_2(\cdot))\in L^2(0,T;\mathbb R^2)$ satisfying $(v_1(t),v_2(t))\in T_{B_{\mathbb R^2}(O_{\mathbb R^2},1)}^\flat((1,0))=\{(v_1,v_2)\in\mathbb R^2|\;v_1\leq0\}$ for almost every $t\in[0,T]$,
where $p^\ell(\cdot)=(p_1^\ell(\cdot),p_2^\ell(\cdot))$ satisfies
\begin{align*}
\left\{\begin{array}{l}
\dot p_1^\ell(t)=\dot p_2^\ell(t)=0,\quad a.e.\,t\in(0,T),
\\p_1^\ell(0)=-\ell_{\psi_1},\;p_2^\ell(0)=-\ell_\phi,\;p_1^\ell(T)=0,\;p_2^\ell(T)=\frac{2\ell_{0_2}T}{1+T^2}+\ell_{\psi_2}.
\end{array}
\right.
\end{align*}
Thus, $\ell=(\ell_{0_1},\ell_{0_2},\ell_\phi,\ell_{\psi_1},\ell_{\psi_2})\in\mathbb R^5\setminus\{O_{\mathbb R^5}\}$ is a Lagrangian multiplier if and only if
\begin{align}\label{elm}
\ell_{0_1},\ell_{0_2},\ell_\phi\in(-\infty,0],\ell_{\psi_1}=0, -\ell_\phi=\frac{2\ell_{0_2}T}{1+T^2}+\ell_{\psi_2}.
\end{align}
By means of Proposition \ref{elsn}, we can check that the function $(\tilde v_1(\cdot), \tilde v_2(\cdot))=(0,1)^\top$ is the singular direction. Applying Theorem \ref{sen} to this direction,
 note that 
\begin{align*}
T_{B_{\mathbb R^2}(O_{\mathbb R^2},1)}^{\flat(2)}\left((1,0), (0,1)  \right)=\left\{(v_1,v_2)\in\mathbb R^2;v_1\leq-\frac{1}{2}\right\}.
\end{align*}
Then, by  Proposition \ref{elsn},  there exists a Lagrangian multiplier $\tilde\ell=(\tilde\ell_{0_1},\tilde\ell_{0_2},\tilde\ell_\phi,\tilde\ell_{\psi_1},\tilde\ell_{\psi_2})\in\mathbb R^5\setminus\{O_{\mathbb R^5}\}$ satisfying  (\ref{elm}) such that
\begin{align}\label{soe}\begin{array}{l}
\displaystyle 2\tilde\ell_\phi\int_0^T\Big(\sigma(t)+\left( 1+\frac{1}{2}K(\bar x_1(t),\bar x_2(t)) \right)X_{\tilde v,1}(t)^2-4X_{\tilde v,1}(t)\Big)dt
\\[2mm]\displaystyle-2\left( \tilde\ell_{0_1}+\frac{1}{1+T^2}\tilde\ell_{0_2} \right)X_{\tilde v,1}(T)^2\leq 0,
\end{array}
\end{align}
holds for all $\sigma(\cdot)\in L^2(0,T;\mathbb R)$ such that $\sigma(t)\leq-\frac{1}{2}$ a.e. $t\in[0,T]$, where $X_{\tilde v,1}(t) =\int_0^t\ln^2(1+\tau^2)d\tau$ for $t\in[0,T]$,
and $K(\bar x_1(t),\bar x_2(t))$ denotes the sectional curvature  of $M$ at $(\bar x_1(t),\bar x_2(t))$, with the expression  
$K(\bar x_1(t),\bar x_2(t))=\frac{4(1-t^4)}{(1+6t^2+t^4)^2}$.

The inequality (\ref{soe}) does not hold, which implies that $\bar u(\cdot)$ is not a weak Pareto optimal solution.  In fact, if $\tilde\ell_\phi < 0$, we can choose $\sigma(t)$ (with $t\in[0,T]$) to be sufficiently small such that (\ref{soe}) does not hold. If $\tilde\ell_\phi = 0$, from (\ref{elm}) we know that $(\tilde\ell_{0_1}, \tilde\ell_{0_2}) \in\left( (-\infty, 0] \times (-\infty, 0] \right)\setminus \{O_{\mathbb R^2}\}$, hence the inequality (\ref{soe}) also does not hold.   $\Box$


\setcounter{equation}{0}

\section{Proofs of the main results }\label{pmr3}
\def\theequation{3.\arabic{equation}}
In this section, we will prove the main results for Problem $(MP)_{r,T}$ and Problem $(MP)_{r}$, respectively.
\subsection{The case of problem $(MP)_{r,T}$}

Since the main ideas of the proofs of Theorem \ref{sen} and Theorem \ref{fon} are the same, I will first provide a detailed proof of Theorem \ref{sen} and then outline the proof of Theorem \ref{fon}.

\subsubsection{Proof of Theorem \ref{sen}}
In this subsection, we assume that all the assumptions in Theorem \ref{sen} hold, and that $(\bar u(\cdot),\bar x(\cdot))$ is a weak Pareto optimal solution of problem $(MP)_{r,T}$.


\medskip

Then, related to the pair $(\bar u(\cdot),\bar x(\cdot))$, we introduce some notations.

Given an index set $I\subset\{0_1,\cdots,0_r,1,\cdots,j\}$, we denote by $\Phi_I=(\bar \phi_{0_1},\cdots,\bar\phi_{0_r},\bar\phi_1,
 \cdots,\bar\phi_j,\psi^\top)^\top$, where 
$\bar\phi_i=\phi_i$ if $i\in I$, and $\bar \phi_i=0$ if $i\notin I$. 

For a singular direction $v(\cdot)$, we denote by $X_v(\cdot)$ the vector field along $\bar x(\cdot)$, satisfying (\ref{sdc}) and (\ref{sdc1}). It follows from \cite[Lemma 4.3]{dc} that $v(\cdot)\in T_{\mathcal U_T\cap L^2(0,T;\mathbb R^m)}^\flat(\bar u(\cdot))$.  For $(W,\sigma(\cdot))\in T_{\bar x(0)}M\times T_{\mathcal U_T\cap L^2(0,T;\mathbb R^m)}^{\flat(2)}(\bar u(\cdot),v(\cdot)))$, denote by $Y^{X_v,W}_{v,\sigma}(\cdot)$ the solution to the following equation
\begin{align}\label{sve12}
\left\{\begin{array}{ll}\nabla_{\dot{\bar x}(t)}Y_{v,\sigma}^{X_v,W}(\mu)=\nabla_xf[t](\mu,Y_{v,\sigma}^{X_v,W}(t))+\nabla_u f[t]\left(\mu,\sigma(t)\right)+\frac{1}{2}\nabla_x^2f[t]\big(\mu,X_v(t),
\\[2mm]\quad\quad\quad\quad\quad\quad\quad\quad X_v(t)\big)+\nabla_u\nabla_x f[t]\left(\mu,X_v(t), v(t)\right)+\frac{1}{2}\nabla_u^2f[t]\left(\mu, v(t),v(t)\right)
\\[2mm]\quad\quad\quad\quad\quad\quad\quad\quad-\frac{1}{2}R(\tilde\mu,X_v(t),f[t],X_v(t)),\;a.e.\,t\in(0,T),\;\forall\,\mu\in T^*M,
\\ [2mm] Y_{v,\sigma}^{X_v,W}(0)=W,
\end{array}
\right.
\end{align}
where we adopt notation $[t]$  introduced by (\ref{abbre}), $\tilde\mu$ denotes  the dual vector of $\mu$ (see (\ref{dcv}) for the definition),   and $\nabla_u^2f[t]$ and $\nabla_u\nabla_xf[t]$ are defined by (\ref{dfusx}).
  Then, $Y^{X_v,W}_{v,\sigma}(\cdot)$ is a vector field along $\bar x(\cdot)$ (i.e., $Y^{X_v,W}_{v,\sigma}(t)\in T_{\bar x(t)}M$ for $t\in[0,T]$).
 Recalling (\ref{fdv}), (\ref{sdv}) and (\ref{msd}),   we set
\begin{align*}
\mathcal K_v=&\Big\{\nabla_1\Phi_{I_A^\prime}(\bar x(0),\bar x(T))(W)+\nabla_2\Phi_{I_A^\prime}(\bar x(0),\bar x(T))(Y_{v,\sigma}^{X_v,W}(T))
\\&+\frac{1}{2}D^2\Phi_{I_A^\prime}(\bar x(0),\bar x(T))(X_v)|\; (W,\sigma(\cdot))\in  T_{\bar x(0)}M\times T_{\mathcal U_T\cap L^2(0,T;\mathbb R^m)}^{\flat(2)}\left(\bar u(\cdot), v(\cdot)  \right)\Big\},
\\\mathcal K^\psi_v=&\Big\{\nabla_1\psi(\bar x(0),\bar x(T))(W)+\nabla_2\psi(\bar x(0),\bar x(T))(Y_{v,\sigma}^{X_v,W}(T))+\frac{1}{2}D^2\psi(\bar x(0),\bar x(T))(X_v)|\;
\\& (W,\sigma(\cdot))\in  T_{\bar x(0)}M\times T_{\mathcal U_T\cap L^2(0,T;\mathbb R^m)}^{\flat(2)}\left(\bar u(\cdot), v(\cdot)  \right)\Big\}.
\end{align*}

\begin{lem}\label{sos}
Assume assumptions $(A1)$ and  $(A2)$ hold, and $(\bar u(\cdot),\bar x(\cdot)) $ is a weak Pareto optimal solution of problem $(MP)_{r,T}$. Let $v(\cdot)$ be a singular direction  with $X_v(\cdot)$ satisfying  (\ref{sdc}) and (\ref{sdc1}). Then,  both $\mathcal K_v$ and $\mathcal K^\psi_v$ are convex. Moreover,  set $Y=(Y_{0_1},\cdots,Y_{0_r},Y_1,\cdots,Y_j)^\top$
with
\begin{align*}
Y_i^v=\left\{\begin{array}{ll}\nabla_1\phi_i(\bar x(0),\bar x(T))(X_v(0))+\nabla_2\phi_i(\bar x(0),\bar x(T))(X_v(T)),&i\in I_A,
\\ 0,&i\notin I_A,
\end{array}
\right.
\end{align*}
and 
$$P_v=(-\infty,\underbrace{0)\times\cdots\times(-\infty,}_{r+j\;\textrm{times}}0)-\Big\{\lambda\Big(\hat\phi_0(\bar x(0),\bar x(T))^\top,\phi(\bar x(0),\bar x(T))^\top)^\top+Y\Big)|\;\lambda> 0\Big\},$$ where $\hat\phi_0(x_1,x_2)\stackrel{\triangle}{=}\phi_0(x_1,x_2)-\phi_0(\bar x(0),\bar x(T))$ for $x_1,x_2\in M$. If there does not exist a hyperplane  separating $\mathcal K_v$ and $P_v\times\{O_{\mathbb R^{k}}\}$ properly, i.e. there does not exist $\ell\in\mathbb R^{r+j+k}\setminus\{O_{\mathbb R^{r+j+k}}\}$ such that
\begin{align}\label{sptak}
\ell^\top\nu
\leq\ell^\top(z^\top,O_{\mathbb R^k})^\top,\quad\forall\,(\nu,z)\in\mathcal K_v \times P_v,
\end{align}
  then there exist $(\sigma_1(\cdot),W^1),\cdots,(\sigma_{k+1}(\cdot),W^{k+1})\in T_{\mathcal U_T\cap L^2(0,T;\mathbb R^m)}^{\flat(2)}\times T_{\bar x(0)}M$ and $\delta_0>0$ such that
\begin{align}
&\label{bcch}\begin{array}{ll}
\ds B_{\mathbb R^k}(O_{\mathbb R^k},\delta_0)
\subset& \textrm{co}\,\{\nabla_1\psi(\bar x(0),\bar x(T))(W^i)+\nabla_2\psi(\bar x(0),\bar x(T))(Y_{v,\sigma_i}^{X_v,W^i}(T))
\\[2mm]&\ds+\frac{1}{2}D^2\psi(\bar x(0),\bar x(T))(X_v)\}_{i=1}^{k+1},
\end{array}
\\& \label{ikp1}\begin{array}{l}
\displaystyle\nabla_1\phi_i(\bar x(0),\bar x(T))(W^\eta)+\nabla_2\phi_i(\bar x(0),\bar x(T))(Y_{v, \sigma_\eta}^{X_v,W^\eta}(T))
\\[2mm]\displaystyle+\frac{1}{2}D^2\phi_i(\bar x(0),\bar x(T))(X_u)< 0,\;i\in I_A^\prime,\;\eta=1,\cdots,k+1,
\end{array}
\end{align}
where $\textrm{co}\,A$ is the convex hull of set $A$, $B_{\mathbb R^k}(O_{\mathbb R^k},\delta_0)$ is the closed ball in $\mathbb R^k$ with center zero $O_{\mathbb R^k}$ and radius $\delta_0$, and $D^2\psi(\bar x(0),\bar x(T))(X_v)$ and $D^2\phi_i(\bar x(0),\bar x(T))(X_v)$ are defined by  (\ref{msd}).
\end{lem}

{\it Proof.}\;
First, we show that both $\mathcal K_v$ and $\mathcal K_v^\psi$ are convex. In fact, we only have to show that  the following set
\begin{align*}
\mathcal R^2(T)=\left\{Y_{v,\sigma}^{X_v,W}(T)|\;\sigma(\cdot)\in T_{\mathcal U_T\cap L^2(0,T;\mathbb R^m)}^{\flat(2)}(\bar u(\cdot),v(\cdot)), \,W\in T_{\bar x(0)}M\right\}
\end{align*}
is convex. 
To this end, take an orthonormal basis for $T_{\bar x(0)}M$ :$\{e_1,\cdots,e_n\}\subset T_{\bar x(0)}M$, i.e. $\langle e_i,e_l\rangle=\delta_i^l$ for $i,l=1,\cdots,n$. Denote by $\{d_1,\cdots,d_n\}\subset T_{\bar x(0)}^*M$ the dual basis of $\{e_1,\cdots,e_n\}$ at $\bar x(0)$, i.e. $d_i(e_l)=\delta_i^l$ for $i,l=1,\cdots,n$. For each $t\in[0,T]$, set $e_i(t)=L_{\bar x(0)\bar x(t)}^{\bar x(\cdot)}e_i$ and $d_i(t)=L_{\bar x(0)\bar x(t)}^{\bar x(\cdot)}d_i$ for $i=1,\cdots,n$, which are the parallel translation of $\{e_1,\cdots,e_n\}$ and $\{d_1,\cdots,d_n\}$ from $\bar x(0)$ to $\bar x(t)$, respectively. It follows from
\cite[(2.6) \& (2.7)]{cdz} that, $\{e_1(t),\cdots,e_n(t)\}$ is an orthonormal basis for $T_{\bar x(t)}M$, and $d_i(t)(e_l(t))=\delta_i^l$ for $i,l=1,\cdots,n$ and $t\in[0,T]$.

Then, we shall express  $\mathcal R^2(T)$ 
by means of $\{e_i(\cdot)\}_{i=1}^n$ and $\{d_i(\cdot)\}_{i=1}^n$. To this end, take any $W\in T_{\bar x(0)}M$ and $\sigma(\cdot)=\left(\sigma^1(\cdot),\cdots,\sigma^m(\cdot)\right)\in L^2(0,T;\mathbb R^m)$ with $\sigma(t)\in T_U^{\flat(2)}(\bar u(t), v(t))$ a.e. $t\in(0,T)$. Set
\begin{align*}
& Y_{v,\sigma,i}^{X_v, W}(t)=Y_{v,\sigma}^{X_v, W}(t)(d_i(t)),\quad i=1,\cdots,n,\;t\in[0,T],
\\ & F_{il}^x(t)=\nabla_xf[t](d_i(t),e_l(t)),\quad i,l=1,\cdots,n,\;t\in[0,T],
\\ & F_{il}^{u}(t)=\nabla_uf[t](d_i(t),E_l),\quad i=1,\cdots,n,\;l=1,\cdots,m,\;t\in[0,T],
\\& S_{v,i}^{X_v}(t)=\frac{1}{2}\nabla_x^2f[t]\left(d_i(t), X_v(t),X_v(t) \right)+\nabla_u\nabla_x f[t]\left(d_i(t), X_v(t),v(t)  \right)
\\&\qquad\quad\quad+\frac{1}{2}\nabla_u^2f[t]\left(d_i(t),v(t),v(t)  \right)-\frac{1}{2} R\left(  e_i(t), X_v(t), f[t],X_v(t) \right),
\end{align*}
where $E_l=(\underbrace{0,\cdots,1}_{l},\cdots,0)^\top\in\mathbb R^m$.
Thus, we obtain from the definition of tensors (see \cite[Section 2.2]{cdz}) and the properties of $\{e_i(\cdot)\}_{i=1}^n$ and $\{d_i(\cdot)\}_{i=1}^n$ that
\begin{align*}
& Y_{v,\sigma}^{X_v, W}(t)=\sum_{i=1}^n Y_{v,\sigma,i}^{X_v, W}(t)e_i(t),\quad \nabla_xf[t]=\sum_{i,l=1}^n F_{il}^x(t)e_i(t)\otimes d_l(t),\;t\in[0,T],
\\ &\nabla_uf[t](\cdot,\sigma(t))=\sum_{i=1}^n\sum_{l=1}^mF_{il}^{u}(t)\sigma^l(t) e_i(t),\quad t\in[0,T].
\end{align*}
Set 
\begin{align*}
&\vec Y_{v,\sigma}^{X_v,W}(t)=\left( Y_{v,\sigma,1}^{X_v,W}(t),\cdots,   Y_{v,\sigma,n}^{X_v,W}(t) \right)^\top,\;F_x(t)=(F_{il}^x(t))_{n\times n},\; F_u(t)=(F_{il}^u(t))_{n\times m},
\\  & \vec S_v^{X_v}(t)=\left( S_{v,1}^{X_v}(t),\cdots,S_{v,n}^{X_v}(t) \right)^\top;\;\vec W=\left(W(d_1),\cdots,W(d_n)\right)^\top,
\end{align*}
for $t\in[0,T]$.  Note that
\begin{align*}
\nabla_{\dot{\bar x}(t)}Y_{v,\sigma}^{X_v,W}=\sum_{i=1}^n\nabla_{\dot{\bar x}(t)}\Big(Y_{v,\sigma,i}^{X_v,W}(\cdot)e_i(\cdot)\Big)=\sum_{i=1}^n\dot Y_{v,\sigma,i}^{X_v,W}(t)e_i(t),\quad\forall\,t\in[0,T],
\end{align*}
where we have used the definition of the parallel translation (see \cite[Section 2.2]{cdz}) and the definition of $\{e_i(\cdot)\}_{i=1}^n$.  In equation (\ref{sve12}), for almost every $t\in(0,T)$,  by   assigning the values of $e_1(t),\cdots,e_n(t)$ to $\mu$ respectively, we can obtain that the vector-valued function $\vec Y_{v,\sigma}^{X_v,W}(\cdot)$ satisfies the following ordinary differential equation:
\begin{align*}
\left\{\begin{array}{l}\dot{\vec Y}_{v,\sigma}^{X_v,W}(t)=F_x(t)\vec Y_{v,\sigma}^{X_v,W}(t)+F_u(t)\sigma(t)+\vec S_v^{X_v}(t),\quad a.e.\,t\in(0,T),
 \\ \vec Y_{v,\sigma}^{X_v,W}(0)=\vec W.
\end{array}
\right.
\end{align*}
Thus, to show $\mathcal R^2(T)$ is convex, it suffices to prove that the set
$\vec{\mathcal R}^2(T)=\{\vec Y_{v,\sigma}^{X_v,W}(T)|\, W\in T_{\bar x(0)} M,\,\sigma(\cdot)\in T_{\mathcal U_T\cap L^2(0,T;\mathbb R^m)  }^{\flat(2)}(\bar u(\cdot),v(\cdot))\}$
is convex.
Indeed, by the definition of  $T_{\mathcal U_T\cap L^2(0,T;\mathbb R^m)}^{\flat(2)}(\bar u(\cdot),v(\cdot))$ and the convexity of $\mathcal U_T$, we conclude that $T_{\mathcal U_T\cap L^2(0,T;\mathbb R^m)}^{\flat(2)}(\bar u(\cdot),v(\cdot))$ is convex.  Consequently, $\vec{\mathcal R}^2(T)$ is convex.


Second, we shall show that zero belongs to the interior of $\mathcal K^\psi_v$. By contradiction, we assume it were not true.   For the case that  the interior of $\mathcal K^\psi_v$ is not empty, since $\mathcal K_v$ and $\mathcal K^\psi_v$ are both convex,  we obtain from \cite[Theorem 11.1,
p.95 $\&$ Theorem 11.3, p.97]{rock} that, there exists $\xi\in\mathbb R^k\setminus\{O_{\mathbb R^k}\}$ such that   $\xi^\top\nu\leq 0$ for all $\nu\in \mathcal K^\psi_v$, which immediately implies that $(O_{\mathbb R^{r+j}},\xi^\top)\tilde\nu\leq (O_{\mathbb R^{r+j}},\xi^\top)(z^\top,O_{\mathbb R^k})^\top$ holds for all $\tilde\nu\in\mathcal K_v$ and $z\in P_v$. Consequently, a contradiction follows. For the case that the interior of $\mathcal K^\psi_v$ is empty,   the dimension of the subspace spanned by $\mathcal K^\psi_v$ is less than $k$, and consequently, there exists  $\hat\xi\in\mathbb R^k\setminus\{O_{\mathbb R^k}\}$ such that $\hat\xi^\top \nu=0$ for all $\nu\in \mathcal K^\psi_v$, which implies $(O_{\mathbb R^{r+j}},\hat\xi^\top)\hat\nu=(O_{\mathbb R^{r+j}},\hat\xi^\top)(z^\top,O_{\mathbb R^k})^\top$ holds for all $\hat\nu\in\mathcal K_v$ and $z\in P_v$. A contradiction follows.

Third, we obtain from the above paragraph that,  
there exists $\hat\delta_0>0$ and
\\
 $\{(\hat\sigma_\eta(\cdot),  \hat W^\eta)\}_{\eta=1}^{k+1}\subset T_{\mathcal{U}\cap L^2(0,T;\mathbb{R}^m)}^{\flat(2)}(\bar{u}(\cdot),v(\cdot)) \times T_{\bar x(0)}M$ such that
(\ref{bcch}) holds with $\delta_0$ and 
$\{(\sigma_\eta(\cdot),  W^\eta)\}_{\eta=1}^{k+1}$ replaced respectively by  $\hat\delta_0$ and $\{(\hat\sigma_\eta(\cdot),   \hat W^\eta)\}_{\eta=1}^{k+1}$. Moreover, we obtain from \cite[Theorem 11.3, p.97]{rock} that $\mathcal K_v \cap (P_v\times \{O_{\mathbb R^k}\})\not=\emptyset$. In other words,  there exists $(\sigma_0(\cdot),  W^0)\in  T_{\mathcal{U}\cap L^2(0,T;\mathbb{R}^m)}^{\flat(2)}(\bar{u}(\cdot),v(\cdot)) \times T_{\bar x(0)}M$ and  $z_i<0$ (with $i=0,1,\cdots,j$) such that
\begin{align}\label{ikp}\begin{array}{l}
\nabla_1\phi_i(\bar x(0),\bar x(T))(W^0)+\nabla_2\phi_i(\bar x(0),\bar x(T))(Y_{v,\sigma_0}^{X_v,W^0}(T))
+\frac{1}{2}D^2\phi_i(\bar x(0),\bar x(T))(X_v)=z_i,\;i\in I_A^\prime,
\\ [2mm]\nabla_1\psi(\bar x(0),\bar x(T))(W^0)+\nabla_2\psi(\bar x(0),\bar x(T))(Y_{v,\sigma_0}^{X_v,W^0}(T))+\frac{1}{2}D^2\psi(\bar x(0),\bar x(T))(X_v)=0.
\end{array}
\end{align}
We take $\theta\in(0,1)$ such that
\begin{align*}
&(1-\theta)\Big(\nabla_1\phi_i(\bar x(0),\bar x(T))(\hat W^\eta)+\nabla_2\phi_i(\bar x(0),\bar x(T))(Y_{v,\hat\sigma_\eta}^{X_v,\hat W^\eta}(T))
\\&+\frac{1}{2}D^2\phi_i(\bar x(0),\bar x(T))(X_v)\Big)
+\theta \Big(\nabla_1\phi_i(\bar x(0),\bar x(T))( W^0)
\\&+\nabla_2\phi_i(\bar x(0),\bar x(T))(Y_{v,\sigma_0}^{X_v, W^0}(T))
+\frac{1}{2}D^2\phi_i(\bar x(0),\bar x(T))(X_v)\Big)<0
\end{align*}
holds for $\eta=1,\cdots,k+1$ and $i\in I_A^\prime$. 
Set $W^\eta=(1-\theta)\hat W^\eta+\theta W^0$ and $\sigma_\eta(\cdot)=\\ (1-\theta)\hat\sigma_\eta(\cdot)+\theta\sigma_0(\cdot)$ for $\eta=1,\cdots,k+1$, and  $\delta_0=(1-\theta)\hat\delta_0$.
Then, 
(\ref{bcch}) and (\ref{ikp1}) follows.  The proof is concluded. $\Box$

\medskip

\begin{proposition}\label{pmt}
Assume all the assumptions in Theorem  \ref{sen} hold   and $\phi_0(\bar x(0),\bar x(T)) =O_{\mathbb R^r}$.  Let $v(\cdot)$ be a singular direction,  with $X_v$ satisfying (\ref{sdc}) and (\ref{sdc1}). Then, there exists $\ell\in\mathbb R^{r+j+k}\setminus\{O_{\mathbb R^{r+j+k}}\}$ satisfying (\ref{sptak}).
\end{proposition}

{\it Proof.}\; The proof is by contradiction. If there did not exist $\ell\in \mathbb R^{r+j+k}\setminus\{O_{\mathbb R^{r+j+k}}\}$ satisfying  (\ref{sptak}), we would obtain from Lemma \ref{sos}  that, there exist $\{(\sigma_\eta(\cdot), W^\eta)\}_{\eta=1}^{k+1}\subset T_{\mathcal U_T\cap L^2(0,T;\mathbb R^m)}^{\flat(2)}(\bar u(\cdot),v(\cdot)) \times T_{\bar x(0)}M$ and $\delta_0>0$ such that (\ref{bcch}) and (\ref{ikp1}) hold.

Let $\{h_\zeta\}_{\zeta\geq1}$ be a sequence such that $h_\zeta\to0^+$ as $\zeta\to+\infty.$ We obtain from the characterization of the second-order adjacent subset that, for each $\eta\in\{1,\cdots,k+$ $1\}$, there eixsts $\{\sigma_{\eta}^\zeta(\cdot)\}_{\zeta=1}^{+\infty}\subset L^2(0,T;\mathbb{R}^m)$ such that $\bar{u}(\cdot)+h_\zeta v(\cdot)+h_\zeta^2\sigma_\eta^\zeta(\cdot)\in\mathscr{U}_T$ for all $\zeta\geq1$ and $\sigma_{\eta}^{\zeta}(\cdot)\to\sigma_{\eta}(\cdot)$ in $L^2(0,T;\mathbb{R}^m)$ as $\zeta\to+\infty$.
For any $\nu=(\nu_1,\cdots,\nu_{k+1})^\top\in\mathbb R^{k+1}$ satisfying 
\begin{align}\label{clc}
\nu_i\geq 0,\quad i=1,\cdots,k+1;
\quad\sum_{i=1}^{k+1}\nu_i=1,
\end{align}
 we define
$
u_{v,\nu}^\zeta(\cdot)=\bar{u}(\cdot)+h_\zeta v(\cdot)+h_\zeta^2\sum_{\eta=1}^{k+1}\nu_\eta\sigma_\eta^\zeta(\cdot)$ for all $\zeta\geq1$. 
Then, $u_{v,\nu}^\zeta(\cdot)\in\mathcal U_T$. Denote by $x_{v,\nu}^\zeta(\cdot)$ the solution to (\ref{25}) with control $u_{v,\nu}^\zeta(\cdot)$ and the initial state $x_{v,\nu}^\zeta(0)=\exp_{\bar x(0)}\left( h_\zeta X_v(0)+h_\zeta^2\sum_{\eta=1}^{k+1}\nu_\eta W^\eta \right)$. It follows from \cite[Proposition 3.1]{dc} that, there exists  $\zeta_2\geq 1$ such that  we can define
$D_\nu^\zeta(t)=exp_{\bar x(t)}^{-1}x_{v,\nu}^\zeta(t)$ for all $t\in[0,T]$ and $\zeta\geq \zeta_2$.
For each $t\in[0,T]$, and $\zeta\geq\zeta_2$, set
$
\alpha_{\nu}^\zeta(\theta;t)=\exp_{\bar x(t)}\left( \theta D_\nu^\zeta(t) \right)$ for $\theta\in[0,1]$. 
By \cite[Lemma 2.1]{cdz}, the curve $\alpha_\nu^\zeta(\cdot;t)$ is the shortest geodesic connecting $\bar x(t)$ and $x_{v,\nu}^\zeta(t)$. We obtain from \cite[Lemma 2.1]{cdz} and the definitions of the geodesic (see \cite[(2.1)]{cdz}) and parallel translation (see \cite[Section 2.2]{cdz}) that
\begin{align}\label{dsg}
\frac{\partial}{\partial\theta}\Big|_0\alpha_\nu^\zeta(\theta;t)= D_\nu^\zeta(t);\quad \frac{\partial}{\partial\theta}\alpha_\nu^\zeta(\theta;t)= L_{\bar x(t)\alpha_\nu^\zeta(\theta;t)} D_\nu^\zeta(t),\quad\theta\in[0,1]
\end{align}
holds for $t\in[0,T]$. For 
 $\varphi(\cdot)=\phi_{0_1}(\cdot),\cdots,
\phi_{0_r}(\cdot),\phi_1(\cdot),\cdots,\phi_j(\cdot),\psi(\cdot)$, by using the Newton-Leibniz formula, interchanging the order of integration and (\ref{dsg}), we can obtain 
\begin{align}
&\varphi\left( x_{v,\nu}^\zeta(0),x_{v,\nu}^\zeta(T)  \right)-\varphi\left(\bar x(0),\bar x(T)  \right)\nonumber
\\&=\int_0^1\frac{\partial}{\partial\theta}\varphi\left( \alpha_\nu^\zeta(\theta;0), \alpha_\nu^\zeta(\theta;T) \right)d\theta\nonumber
\\&=h_\zeta\Big(\nabla_1\varphi\left(  \bar x(0),\bar x(T)\right)(X_v(0))+\nabla_2\varphi\left(  \bar x(0),\bar x(T)\right)(X_v(T))\Big)\nonumber
\\&\quad +\int_0^1\Bigg(\nabla_1\varphi \left( \alpha_\nu^\zeta(\theta;0), \alpha_\nu^\zeta(\theta;T) \right)\left( \frac{\partial}{\partial\theta}\alpha_\nu^\zeta(\theta;0) \right)-\nabla_1\varphi(\bar x(0),\bar x(T))\left( \frac{\partial}{\partial\theta}\Big|_{0}\alpha_\nu^\zeta(\theta;0) \right)\nonumber
\\&\quad+\nabla_2\varphi \left( \alpha_\nu^\zeta(\theta;0), \alpha_\nu^\zeta(\theta;T) \right)\left( \frac{\partial}{\partial\theta}\alpha_\nu^\zeta(\theta;T) \right)-\nabla_2\varphi(\bar x(0),\bar x(T))\left( \frac{\partial}{\partial\theta}\Big|_{0}\alpha_\nu^\zeta(\theta;T) \right)\Bigg)d\theta\nonumber
\\&\quad+\nabla_1\varphi(\bar x(0),\bar x(T))(D_\nu^\zeta(0)-h_\zeta X_v(0))+\nabla_2\varphi(\bar x(0),\bar x(T))(D_\nu^\zeta(T)-h_\zeta X_v(T))\nonumber
\\&=h_\zeta\Big(\nabla_1\varphi\left(  \bar x(0),\bar x(T)\right)(X_v(0))+\nabla_2\varphi\left(  \bar x(0),\bar x(T)\right)(X_v(T))\Big)+h_\zeta^2\sum_{\eta=1}^{k+1}\nu_\eta P_{\varphi,\eta}\nonumber
\\&\quad+h_\zeta^2  \sum_{\eta=1}^{k+1}\nu_\eta  \nabla_2\varphi(\bar x(0),\bar x(T))\left(Y_{v,\sigma_\eta^\zeta}^{X_v,W^\eta}(T)- Y_{v,\sigma_\eta}^{X_v, W^\eta}(T)\right)\nonumber
\\&\quad+\nabla_1\varphi(\bar x(0),\bar x(T))\left( D_\nu^\zeta(0)-h_\zeta X_v(0)-h_\zeta^2\sum_{\eta=1}^{k+1}\nu_\eta W^\eta  \right)\nonumber
\\&\quad+\nabla_2\varphi(\bar x(0),\bar x(T))\left( D_\nu^\zeta(T)-h_\zeta X_v(T)-h_\zeta^2 Y_{v,\sum_{\eta=1}^{k+1}\nu_\eta\sigma_\eta^\zeta}^{X_v,\sum_{\eta=1}^{k+1}\nu_\eta W^\eta}(T)\right)\nonumber
\\&\quad+\int_0^1\int_0^\theta\frac{\partial}{\partial\tau}\left(\nabla_1\varphi\left(\alpha_\nu^\zeta(\tau;0),  \alpha_\nu^\zeta(\tau;T)\right)\frac{\partial}{\partial\tau}\alpha_\nu^\zeta(\tau;0)\right.\nonumber
\\&\quad\left.+\nabla_2\varphi\left(\alpha_\nu^\zeta(\tau;0),  \alpha_\nu^\zeta(\tau;T)\right)\frac{\partial}{\partial\tau}\alpha_\nu^\zeta(\tau;T)\right)d\tau d\theta-\frac{1}{2}D^2\varphi(\bar x(0),\bar x(T))(X_v)\nonumber
\\\label{sovv} &=h_\zeta\Big(\nabla_1\varphi\left(  \bar x(0),\bar x(T)\right)(X_v(0))+\nabla_2\varphi\left(  \bar x(0),\bar x(T)\right)(X_v(T))\Big)
+h_\zeta^2\sum_{\eta=1}^{k+1}\nu_\eta P_{\varphi,
\eta}+R_{v,\nu}^\zeta,
\end{align}
where
\begin{align*}
P_{\varphi,
\eta}= \nabla_1\varphi(\bar x(0),\bar x(T))(W^\eta)+\nabla_2\varphi(\bar x(0),\bar x(T)\left(  Y_{v,\sigma_\eta}^{X_v, W^\eta}(T) \right)+\frac{1}{2}D^2\varphi(\bar x(0),\bar x(T))(X_v),
\end{align*}
and
\begin{align*}
R_{v,\nu}^{\varphi,\zeta}
&= h_\zeta^2  \sum_{\eta=1}^{k+1}\nu_\eta  \nabla_2\varphi(\bar x(0),\bar x(T))\left(Y_{v,\sigma_\eta^\zeta}^{X_v,W^\eta}(T)- Y_{v,\sigma_\eta}^{X_v, W^\eta}(T)\right)
\\&\quad +\nabla_1\varphi(\bar x(0),\bar x(T))\left( D_\nu^\zeta(0)-h_\zeta X_v(0)-h_\zeta^2\sum_{\eta=1}^{k+1}\nu_\eta W^\eta \right)
\\&\quad+\nabla_2\varphi(\bar x(0),\bar x(T))\left( D_\nu^\zeta(T)-h_\zeta X_v(T)-h_\zeta^2\sum_{\eta=1}^{k+1}\nu_\eta Y_{v,\sigma_\eta^\zeta}^{X_v,W^\eta}(T) \right)
\\&\quad+\int_0^1(1-\tau)\Bigg( \nabla_1^2\varphi\big( \alpha_\nu^\zeta(\tau;0),\alpha_{\nu}^\zeta(\tau;T)  \big)  \Big( L_{\bar x(0)\alpha_\nu^\zeta(\tau;0)}D_\nu^\zeta(0),
\\&\quad L_{\bar x(0)\alpha_\nu^\zeta(\tau;0)}D_\nu^\zeta(0)  \Big)+2\nabla_2\nabla_1\varphi\big( \alpha_\nu^\zeta(\tau;0),\alpha_{\nu}^\zeta(\tau;T)  \big) \Big( L_{\bar x(0)\alpha_\nu^\zeta(\tau;0)}D_\nu^\zeta(0),
\\&\quad L_{\bar x(T)\alpha_\nu^\zeta(\tau;T)}D_\nu^\zeta(T)  \Big)+\nabla_2^2\varphi\big( \alpha_\nu^\zeta(\tau;0),\alpha_{\nu}^\zeta(\tau;T)  \big)\Big(L_{\bar x(T)\alpha_\nu^\zeta(\tau;T)}D_\nu^\zeta(T),
\\&\quad L_{\bar x(T)\alpha_\nu^\zeta(\tau;T)}D_\nu^\zeta(T)\Big)-\nabla_1^2\varphi(\bar x(0),\bar x(T))\big(D_\nu^\zeta(0),D_\nu^\zeta(0)   \big)
\\&\quad-2\nabla_2\nabla_1\varphi(\bar x(0),\bar x(T))\big(D_\nu^\zeta(0),D_\nu^\zeta(T)   \big)
-\nabla_2^2\varphi(\bar x(0),\bar x(T))\big(D_\nu^\zeta(T),D_\nu^\zeta(T)   \big)\Bigg)d\tau
\\&\quad+\frac{1}{2}\Big(  \nabla_1^2\varphi(\bar x(0),\bar x(T))\big(D_\nu^\zeta(0),D_\nu^\zeta(0)\big)
+2\nabla_2\nabla_1\varphi(\bar x(0),\bar x(T))\big(D_\nu^\zeta(0),D_\nu^\zeta(T)   \big)
\\&\quad+\nabla_2^2\varphi(\bar x(0),\bar x(T))\big(D_\nu^\zeta(T),D_\nu^\zeta(T)   \big)-\frac{1}{2}h_\zeta^2D^2\varphi(\bar x(0),\bar x(T))(X_v).
\end{align*}
Recalling (\ref{sve12}), we obtain from direct  computations that
\begin{align*}
\left\{\begin{array}{l}\nabla_{\dot{\bar x}(t)}\left( Y_{v,\sigma_\eta^\zeta}^{X_v, W^\eta}- Y_{v,\sigma_\eta}^{X_v, W^\eta}\right)(\mu)=\nabla_xf[t]\left(\mu,Y_{v,\sigma_\eta^\zeta}^{X_v, W^\eta}(t)-Y_{v,\sigma_\eta}^{X_v, W^\eta}(t)\right)
\\\quad\quad\quad\quad\quad\quad\quad\quad\quad\quad\quad\quad\quad\quad+\nabla_uf[t]\left(\mu,\sigma_\eta^\zeta(t)-\sigma_\eta(t)\right),\;a.e.\,t\in(0,T),\,\mu\in T^*M,
\\ Y_{v,\sigma_\eta^\zeta}^{X_v, W^\eta}(0)-Y_{v,\sigma_\eta}^{X_v, W^\eta}(0)=0.
\end{array}
\right.
\end{align*}
In particular, for almost every $t\in(0,T)$, we  take $\mu=\mu(t)\stackrel{\triangle}{=}\widetilde{Y_{v,\sigma_\eta^\zeta}^{X_v, W^\eta}- Y_{v,\sigma_\eta}^{X_v, W^\eta}} $  in the first line of the above equation, which is the dual covector of $Y_{v,\sigma_\eta^\zeta}^{X_v, W^\eta}(t)- Y_{v,\sigma_\eta}^{X_v, W^\eta}(t)$. Then, we integrate it  
  over  $[0,s]$ (with $s\in[0,T]$) and obtain from (\ref{dcv}) and  \cite[Theorem 3.6, p.55;\, Corollary 3.3, p.54]{c} that
\begin{align*}
&\left| Y_{v,\sigma_\eta^\zeta}^{X_v, W^\eta}(s)- Y_{v,\sigma_\eta}^{X_v, W^\eta}(s)\right|^2
\\&=\int_0^s\frac{d}{dt}\left| Y_{v,\sigma_\eta^\zeta}^{X_v, W^\eta}(t)- Y_{v,\sigma_\eta}^{X_v, W^\eta}(t)\right|^2dt
\\&=2\int_0^s\nabla_{\dot{\bar x}(t)}\left(Y_{v,\sigma_\eta^\zeta}^{X_v, W^\eta}- Y_{v,\sigma_\eta}^{X_v, W^\eta}\right)(\mu(t))dt
\\&=2\int_0^s\left(\nabla_xf[t]\left(\mu(t),Y_{v,\sigma_\eta^\zeta}^{X_v, W^\eta}(t)-Y_{v,\sigma_\eta}^{X_v, W^\eta}(t)\right)
+\nabla_uf[t]\left(\mu(t),\sigma_\eta^\zeta(t)-\sigma_\eta(t)\right)\right)dt.
\end{align*}
By \cite[Lemma 4.1]{cdz} and assumption $(A1)$,  there exists $K_1>0$ such that the following    inequality
$$\left| Y_{v,\sigma_\eta^\zeta}^{X_v, W^\eta}(s)- Y_{v,\sigma_\eta}^{X_v, W^\eta}(s)\right|^2
\leq 3K\int_0^s 
\left| Y_{v,\sigma_\eta^\zeta}^{X_v, W^\eta}(t)- Y_{v,\sigma_\eta}^{X_v, W^\eta}(t)\right|^2dt+K_1\|\sigma_\eta^\zeta-\sigma_\eta\|_{L^2(0,T;\mathbb R^m)}^2$$  holds  for $s\in[0,T]$. Applying Gronwall's inequality to the above inequality, we obtain that there exists $K_2>0$ such that
\begin{align*}
\left| Y_{v,\sigma_\eta^\zeta}^{X_v, W^\eta}(s)- Y_{v,\sigma_\eta}^{X_v, W^\eta}(s)\right|^2\leq K_2\|\sigma_\eta^\zeta-\sigma_\eta\|_{L^2(0,T;\mathbb R^m)}^2,\quad s\in[0,T].
\end{align*}
According to \cite[Proposition 3.1]{dc}, Lebesgue's convergence theorem, (\ref{ikp1}), (\ref{sdc1}),  (\ref{sovv}) and  the above inequality, we can find $\zeta_3\geq\zeta_2$, which does not depend on $\nu$,  such that
\begin{align}\label{ada}
|R_{v,\nu}^{\psi,\zeta}|\leq\delta_0 h_\zeta^2\;\textrm{and}\;\phi_i\left( x_{v,\nu}^\zeta(0), x_{v,\nu}^\zeta(T) \right)<0\;\textrm{for}\,i=0_1,\cdots,0_r,1,\cdots,j,\;\zeta\geq\zeta_3,
\end{align}
where $\delta_0$ is given by (\ref{bcch}). So, by using (\ref{bcch}), (\ref{ada}) and (\ref{sdc1}),   we define a continuous map $\mathcal P:\textrm{co}\,\{P_{\psi,\eta}|\,\eta=1,\cdots,k+1\}\to \textrm{co}\,\{P_{\psi,\eta}|\,\eta=1,\cdots,k+1\}$ as 
$
\mathcal P\left(  \sum_{\eta=1}^{k+1}\nu_\eta P_{\psi,\eta}\right)=-\frac{\psi\left( x_{v,\nu}^{\zeta_3}(0),  x_{v,\nu}^{\zeta_3}(T)\right)}{h_{\zeta_3}^2}+\sum_{\eta=1}^{k+1}\nu_\eta P_{\psi,\eta}$, 
where $\nu=(\nu_1,\cdots,\nu_{k+1})^\top\in\mathbb R^{k+1}$ satisfies  (\ref{clc}). According to Brouwer's fixed  point theorem, there exists $\nu^0=(\nu_1^0,\cdots,\nu_{k+1}^0)^\top\in\mathbb R^{k+1}$ satisfying (\ref{clc}) such that
$
\mathcal P\left(   \sum_{\eta=1}^{k+1}\nu^0_\eta P_{\psi,\eta}\right)=\sum_{\eta=1}^{k+1}\nu^0_\eta P_{\psi,\eta}$, 
which implies
$\psi\left( x_{v,\nu^0}^{\zeta_3}(0),  x_{v,\nu^0}^{\zeta_3}(T)\right) =0$. And according to (\ref{ada}), we get \\
$
\phi_i\left( x_{v,\nu^0}^{\zeta_3}(0),  x_{v,\nu^0}^{\zeta_3}(T)\right)<0$ for $ i=0_1,\cdots,0_r,1,\cdots,j$. 
Thus, $u_{v,\nu^0}^{\zeta_3}(\cdot)$ is an admissible control and the corresponding value $\phi_{0_i}\left( x_{v,\nu^0}^{\zeta_3}(0),  x_{v,\nu^0}^{\zeta_3}(T)\right) <0$ for $i=1,\cdots,r$, which contradicts the assumption that $(\bar u(\cdot),\bar x(\cdot))$ is a weak Pareto optimal solution. The proof is completed. $\Box$

\medskip

{\it Proof of Theorem \ref{sen}.} First, let
$\hat\phi_0(\cdot,\cdot)$ be defined in Lemma \ref{sos}.  Then, $(\bar u(\cdot),\bar x(\cdot))$ is also a weak Pareto optimal solution for the following problem:
\begin{description}
\item[$(\widehat{\textrm{MP}})$] Minimize $\hat\phi_0(x(0),x(T))$ subject to the state equation (\ref
{25}), the control constraint $u(\cdot)\in\mathcal U_T$ and endpoints-constraints (\ref{n22}).
\end{description}
It is obvious that  $\hat\phi_0(\bar x(0),\bar x(T))=O_{\mathbb R^r}$. We obtain from  Proposition \ref{pmt} that,  there exists $\ell=(\ell_{0_1},\cdots,\ell_{0_r},\ell_1,\cdots,\ell_j,\ell_\psi^\top)^\top\in\mathbb R^{r+j+k}\setminus\{O_{\mathbb R^{r+j+k}}\}$ such that (\ref{sptak}) holds. We obtain from this inequality that $\ell$   satisfies  (\ref{lm2}) and (\ref{slm}). Recalling the Lagrangian function $\mathcal L_\ell$ defined by (\ref{lfc}), we denote by $p^\ell(\cdot)$  the solution to  (\ref{de}).

Second, we show that inequality (\ref{sonc2}) holds for all $(W,\sigma(\cdot))\in T_{\bar x(0)} M\times \mathcal S_T$.
To this end, we take
 any $(W,\sigma(\cdot))\in T_{\bar x(0)} M\times \mathcal S_T$. It follows from  \cite[Lemma 4.3]{dc} that $\sigma(\cdot)\in T_{\mathcal U_T\cap L^2(0,T;\mathbb R^m)}^{\flat(2)}(\bar u(\cdot),v(\cdot))$.  Letting $z\in P_v$ tend to zero in inequality (\ref{sptak}), and then using integration by parts over $[0,T]$,     we can obtain   from (\ref{de}) and
(\ref{sve12}) that
\begin{align*}
0&\geq-p^\ell(0)\left(Y_{v,\sigma}^{X_v,W}(0)  \right)+p^\ell(T)\left(Y_{v,\sigma}^{X_v,W}(T)  \right)+\frac{1}{2}D^2\mathcal L_\ell(\bar x(0),\bar x(T))(X_v)
\\&=\int_0^T\left(\nabla_{\dot{\bar x}(t)}p^\ell\left( Y_{v,\sigma}^{X_v,W}(t) \right) +p^\ell(t)\left( \nabla_{\dot{\bar x}(t)}Y_{v,\sigma}^{X_v,W}  \right)  \right)dt+\frac{1}{2}D^2\mathcal L_\ell(\bar x(0),\bar x(T))(X_v)
\\&=\int_0^T\Big(\nabla_u H\{t;\ell\}\left(\sigma(t)\right)+\frac{1}{2}\nabla_x^2 H\{t;\ell\}\left( X_v(t),X_v(t)  \right)+\nabla_u\nabla_x H\{t;\ell\}\left( X_v(t),v(t) \right)
\\&\quad+\frac{1}{2}\nabla_u^2 H\{t;\ell\}(v(t),v(t))-\frac{1}{2}R\left( \tilde p^\ell(t),X_v(t),f[t], X_v(t)  \right)\Big)dt
+\frac{1}{2}D^2\mathcal L_\ell(\bar x(0),\bar x(T))(X_v),
\end{align*}
which implies that (\ref{sonc2})  holds.

Third, we claim that  
$
\int_0^T\nabla_u H\{t;\ell\}(u(t))dt\leq 0$ holds for all $u(\cdot)\in L^2(0,T;\mathbb R^m)$ with $u(t)\in T_U^\flat(\bar u(t))$ a.e. $t\in[0,T]$.
 By contradiction, we assume it were not true. Then, there would exist $\tilde v(\cdot)\in L^2(0,T;\mathbb R^m)$ with $\tilde v(t)\in T_U^\flat(\bar u(t))$ a.e. $t\in[0,T]$, such that $\int_0^T\nabla_u H\{t;\ell\}(\tilde v(t))dt>0$. We obtain from \cite[Lemma 4.3]{dc}, \cite[Lemma 2.4]{fdt} and \cite[Proposition 4.2.1, p.138]{AF} that 
 $T_{\mathcal U_T\cap L^2(0,T;\mathbb R^m)}^{\flat(2)}(\bar u(\cdot), v(\cdot))=T_{\mathcal U_T\cap L^2(0,T;\mathbb R^m)}^{\flat(2)}(\bar u(\cdot), v(\cdot))+T_{\mathcal U_T\cap L^2(0,T;\mathbb R^m)}^\flat(\bar u(\cdot))$.
Take a $\sigma(\cdot)\in\mathcal S_T$. Thus, for any $\lambda>0$,  we have $\sigma(\cdot)+\lambda \tilde v(\cdot)\in T_{\mathcal U_T\cap L^2(0,T;\mathbb R^m)}^{\flat(2)}(\bar u(\cdot), v(\cdot))$.
Since $\int_0^T\nabla_u H\{t;\ell\}(\tilde v(t))dt>0$, we can choose $\lambda>0$ large enough such that
\begin{align*}
&\int_0^T\left(\nabla_uH\{t;\ell\} (\sigma(t)+\lambda\tilde v(t))+\frac{1}{2}\nabla_x^2H\{t;\ell\} (X_v(t),X_v(t))+\nabla_u\nabla_xH\{t;\ell\} (X_v(t),v(t))\right. 
\\&
+\frac12\nabla_u^2H\{t;\ell\}(v(t),v(t))-\frac12R\left(\tilde{p}^\ell(t),X_v(t),f[t],X_v(t)\right)\Big)dt
+\frac12D^2\mathcal{L}_\ell(\bar x(0),\bar x(T))(X_v)>0,
\end{align*}
a contradiction follows.

Finally, we show that $\ell$ is a Lagrangian multiplier for problem $(MP)_{r,T}$. Since $\ell$ satisfies 
(\ref{lm2}), we only have to show  
$\nabla_u H\{t;\ell\}(u(t))\leq 0$ a.e. $t\in[0,T]$ for all $u(\cdot)\in L^2(0,T;\mathbb R^m)$ such that $u(t)\in T_U^\flat(\bar u(t))$ a.e. $t\in[0,T]$.
To this end, we take any $u(\cdot)\in L^2(0,T;\mathbb R^m)$ satisfying $u(t)\in T_U^\flat(\bar u(t))$ a.e. $t\in[0,T]$. Let $t_0\in(0,T]$ be a Lebesgue point of $\nabla_u H\{\cdot;\ell\}(u(\cdot))$, and set $u_\epsilon(\cdot)=I_{[t_0-\epsilon,t_0]}(\cdot)u(\cdot)$ for $\epsilon>0$ small enough, where $I_{[t_0-\epsilon,t_0]}(\cdot)$ is the indicator function of the set $[t_0-\epsilon,t_0]$.
Then, $u_\epsilon(\cdot)\in L^2(0,T;\mathbb R^m)$ and satisfies  $u_\epsilon(t)\in T_U^\flat(\bar u(t))$ a.e $t\in[0,T]$, and consequently we have $ \int_0^T\nabla_u H\{t;\ell\}\frac{u_\epsilon(t)}{\epsilon}dt =\frac{1}{\epsilon}\int_{t_0-\epsilon}^{t_0}\nabla_u H\{t;\ell\}(u(t))dt\leq 0$.
Letting $\epsilon$ tend to zero, we obtain  $\nabla_u H\{t_0;\ell\}(u(t_0)) \leq 0$. This completes the proof.  $\Box$

\subsubsection{Outline of the proof of Theorem \ref{fon}}
Following the approach in Theorem \ref{sen}, we divide the argument into three steps.

{\it Step 1.}\; Set $\mathcal K_1\stackrel{\triangle}{=}\{\nabla_1\Phi_{I_A}(\bar x(0),\bar x(T))(V)+\nabla_2\Phi_{I_A}(\bar x(0),\bar x(T))(X_v^V(T))|
\;(v(\cdot),V)\in \mathcal V(\bar u(\cdot))\times T_{\bar x(0)}M\}$ and $ \mathcal K_1^\psi\stackrel{\triangle}{=}\{\nabla_1\psi(\bar x(0),\bar x(T))(V)+\nabla_2\psi(\bar x(0),\bar x(T))(X_v^V(T))|
(v(\cdot),V)\in \mathcal V(\bar u(\cdot))\times T_{\bar x(0)}M\}$,
where $\mathcal V(\bar u(\cdot))=\{v(\cdot)\in L^2(0,T;\mathbb R^m)|\,v(t)\in T_U^\flat(\bar u(t))\;a.e.\,t\in(0,T)\}$, and $X_v^V(\cdot)$ (with $(v(\cdot),V)\in \mathcal V(\bar u(\cdot))\times T_{\bar x(0)}M$ ) is the solution to (\ref{sdc}),
with initial condition $X_v^V(0)=V$. One can show that $\mathcal K_1\subset \mathbb R^{r+j+k}$ and $\mathcal K_1^\psi\subset \mathbb R^{k}$ are both convex.

{\it Step 2.}\;  Set
$$\begin{array}{ll}\mathcal P=&\underbrace{(-\infty,0)\times\cdots\times(-\infty,0)}_{r+j\,\textrm{times}}
-\{\lambda\big(O_{\mathbb R^r}^\top,\phi_1(\bar x(0),\bar x(T)),\cdots,\phi_j(\bar x(0),\bar x(T))\big)\in\mathbb R^{r+j}|\;\lambda\geq 0\}.
\end{array}$$
Similar to the proof of Proposition \ref{pmt}, we can use Brouwer's fixed point theorem and the contradiction argument to show that,
there  exists  $\ell\in\mathbb R^{r+j+k}\setminus\{O_{\mathbb R^{r+j+k}}\}$ such that
\begin{align}\label{foi5}
\ell^\top \nu\leq \ell^\top (z^\top,O_{\mathbb R^k})^\top,\quad\forall\,(\nu,z)\in \mathcal K_1\times \mathcal P,
\end{align}
which immediately implies (\ref{lm2}).

{\it Step 3.}\; 
Let $\mathcal L_\ell$ be the Lagrangian function  defined by (\ref{lfc}), and let $p^\ell(\cdot)$  be the solution to  (\ref{de}). For any $(V,v(\cdot))\in T_{\bar x(0)}M\times \mathcal V(\bar u(\cdot))$, we obtain from \cite[Lemma 4.3]{dc}  that $v(\cdot)\in T_{\mathcal U_T\cap L^2(0,T;\mathbb R^m)}$.  Let $X_V^v(\cdot)$ denote the solution to  equation (\ref{sdc}) with the initial condition $X_V^v(0)=V$. Letting $z\in\mathcal P$ tend to zero in (\ref{foi5}), we use integration by parts over $[0,T]$. From (\ref{de}) and (\ref{sdc}), we obtain $\int_0^T \nabla_uH\{t;\ell\}(v(t))dt\leq 0$. By following the last step of the proof of Theorem \ref{sen}, we obtain (\ref{lm3}). The proof is complete. $\Box$
\subsection{The case of problem $(MP)_r$}

Before proving Theorem \ref{fnct}, we introduce a related optimal control problem. Take two positive numbers $a$ and $A$ such that $a<A$. Consider the following multi-objective optimal control problem:
\begin{description}
\item[$(MP)_{[0,1]}$]  Minimize $\phi_0(y(0), y(1))$ subject to the state equation
\begin{align}\label{se01}
\left\{\begin{array}{ll}
\dot\tau(s)=v(s),& a.e.\,s\in(0,1),
\\\dot y(s)=v(s)f\left( \tau(s),y(s),w(s) \right),& a.e. \,s\in(0,1),
\\
\tau(0)=0,&
\end{array}
\right.
\end{align}
with the control constraint:
$
\{(v,w):[0,1]\to[a,A]\times U|\,(v(\cdot), w(\cdot))\;\textrm{is measurable}\}$, 
and endpoints-constraints:
$
\phi_i(y(0),y(1))\leq0$ for $i=1,\cdots,j$, $\psi(y(0),y(1))=O_{\mathbb R^k}$ and $\tau(0)=0$.
\end{description}
 The corresponding admissible solution and weak Pareto optimal solution are defined by Definition \ref{adoT}.

\begin{lem}\label{e01p}Assume $(\bar u(\cdot), \bar x(\cdot), \bar T)$ is a weak Pareto optimal solution of problem $(MP)_r$, with $\bar T>0$. Set $\bar v(s)\equiv \bar T$, $\bar y(s)=\bar x(\bar T s)$, $\bar \tau(s)=\bar T s$ and $\bar w(s)=\bar u(\bar Ts)$ for $s\in[0,1]$. Let $a<\bar T$ and $A$ be big enough. Then, $(\bar v(\cdot), \bar w(\cdot), \bar \tau(\cdot),\bar y(\cdot))$ is a weak Pareto optimal solution for the problem $(MP)_{[0,1]}$.

\end{lem}

{\it  Proof.}\quad It is obvious that $(\bar v(\cdot),\bar w(\cdot), \bar \tau(\cdot),\bar y(\cdot))$ is an admissible solution of problem $(MP)_{[0,1]}$. By contradiction, assume it is not a weak Pareto optimal solution of problem $(MP)_{[0,1]}$,  then there exists an admissible solution $(v(\cdot),w(\cdot),\tau(\cdot),y(\cdot))$ of problem $(MP)_{[0,1]}$ such that $\phi_{0_i}(y(0),y(1))<\phi_{0_i}(\bar y(0),\bar y(1))$ for $i=1,\cdots,r$. Since $\dot\tau(s)=v(s)\geq a$ a.e. $s\in[0,1]$, $v(\cdot)$ is strictly monotonically increasing. Thus, the function $\tau(\cdot)$ has an inverse function $\xi:[0,\tau(1)]\to[0,1]$ such that $\tau\left(  \xi(t) \right)=t$ for $t\in[0,\tau(1)]$, $\xi\left(  \tau(s) \right)=s$ for  $s\in[0,1]$, $\xi(0)=0$,
and 
$\xi(\cdot)$ is also  strictly monotonically increasing.
Set  $T=\int_0^1v(s)ds=\tau(1)$, and  $x(t)=y(\xi(t))$ and $u(t)=w(\xi(t))$ for $t\in[0,T]$.
We can show that $u:[0,T]\to U$ is measurable, $(u(\cdot), x(\cdot), T) $ is an admissible solution of problem $(MP)_r$, and $\phi_0\left(y(0),y(1)\right)=\phi_0(x(0),x(T))$. Thus, we have $\phi_{0_i}(x(0),x(T))<\phi_{0_i}(\bar x(0),\bar x(\bar T))$ for $i=1,\cdots,r$, which contradicts the condition that $(\bar u(\cdot),\bar x(\cdot),\bar T)$ is a weak Pareto optimal solution of problem $(MP)_r$. The proof is concluded. $\Box$ 

\medskip

Then, we prove Theorem \ref{fnct}.

\medskip

{\it Proof of Theorem \ref{fnct}}. We adopt the notations introduced in Lemma \ref{e01p}. According to  Lemma \ref{e01p}, $(\bar v(\cdot),\bar w(\cdot), \bar \tau(\cdot),\bar y(\cdot))$ is a weak Pareto optimal solution of problem $(MP)_{[0,1]}$, with $(\bar v(\cdot), \bar w(\cdot))\in L^2(0,1;\mathbb R^{m+1})$.  Let $\mathcal L_\ell$ and $H$ are respectively the Lagrangian function defined by  (\ref{lfc}), and   the Hamiltonian function given by (\ref{hcm}). By means of Theorem  \ref{fon}, there exists a  vector  $\ell$ (given by (\ref{lm1}))
and $\ell_\tau\in\mathbb R$
satisfying (\ref{aim}) and $(\ell,\ell_\tau)\not=O_{\mathbb R^{r+j+k+1}}$, such that $H\left(  \bar T s,\bar y(s), p_1^\ell(s),\bar w(s) \right)+p_{0,1}(s)=0$ and $\nabla_uH\left(  \bar T s,\bar y(s), p_1^\ell(s),\bar w(s) \right)(v(s))\leq 0$ hold for almost every $s\in[0,1]$ and all $v(\cdot)\in L^2(0,1;\mathbb R^m)$ with $v(s)\in T_U^\flat(\bar w(s))$ a.e. $s\in[0,1]$,
  where $p_{0,1}(\cdot)$ and $p_1^\ell(\cdot)$ satisfy
\begin{align}\label{det}
\left\{\begin{array}{l}
\displaystyle\dot p_{0,1}(s)=-\bar T\nabla_t H\left(  \bar T s,\bar y(s), p_1^\ell(s),\bar w(s) \right),\quad a.e.\,s\in[0,1],
\\[2mm] \nabla_{\dot{\bar y}(s)}p_1^\ell=-\bar T\nabla_xf\left( \bar T s,\bar y(s),\bar w(s) \right)\left(   p_1^\ell(s),\cdot \right),\quad a.e.\,s\in[0,1],
\\[2mm] p_{0,1}(0)=-\ell_\tau,\;p_{0,1}(1)=0,\;p_1^\ell(0)=-d_1\mathcal L_\ell\left(\bar y(0),\bar y(1)\right),\;p_1^\ell(1)=d_2\mathcal L_\ell\left(\bar y(0),\bar y(1)\right).
\end{array}
\right.
\end{align}
 By setting $t=\bar T s$ and $p^\ell(t)=p_1^\ell(\frac{t}{\bar T} )$ for $t\in[0,\bar T]$, the above conclusion can be equivalently stated as: there exists a nonzero vector $\ell$ satisfying (\ref{aim}), such that  (\ref{font}) 
holds for all $v(\cdot)\in L^2(0,\bar T;\mathbb R^m)$ with $v(t)\in T_U^\flat(\bar u(t))$ a.e. $t\in[0,\bar T]$, where $p^\ell(\cdot)$ satisfies (\ref{detd}). Proof completed.  $\Box$

\medskip

{\it Proof of Theorem \ref{snt}.}\; We adopt all the notations introduced in Lemma \ref{e01p} and the proof of Theorem  \ref{fnct}.
Take any singular direction $(\xi(\cdot), v(\cdot))\in L^2(0,\bar T;\mathbb R\times\mathbb R^m)$. Let $X_{\xi,v}(\cdot)$ satisfies (\ref{foet}) and (\ref{sdt}). Set
\begin{align}\label{vtt}
v_0(s)=\xi(\bar T s),\;w(s)=v(\bar T s),\; y_{v_0}(s)=\int_0^s v_0(\mu)d\mu,\; Y_w(s)=X_{\xi,v}(\bar T s),\; s\in[0,1].
\end{align}
We can verify that $\left(v_0(\cdot), w(\cdot), y_{v_0}(\cdot), Y_w(\cdot)  \right)$ fulfills 
 the following equation
\begin{align*}
\left\{\begin{array}{l}
\dot y_{v_0}(s)=v_0(s),\quad a.e.\,s\in(0,1),
\\\nabla_{\dot{\bar y}(s)}Y_w=\bar T\nabla_xf\left(\bar T s, \bar y(s), \bar w(s)  \right)\left(\cdot, Y_w(s)  \right)+\bar T f_t(\bar Ts,\bar y(s),\bar w(s))y_{v_0}(s)
\\\quad\quad\quad\quad\quad+v_0(s)f\left(\bar T s,\bar y(s), \bar w(s)\right)+\bar T\nabla_u f\left(\bar T s,\bar y(s), \bar w(s)\right)(w(s)),\quad a.e.\,s\in(0,1),
\\ y_{v_0}(0)=0,\quad \nabla_1\psi\left(  \bar y(0),\bar y(1) \right) Y_w(0)+\nabla_2\psi\left(  \bar y(0),\bar y(1) \right) Y_w(1)=0,
\\\nabla_1\phi_i\left(  \bar y(0),\bar y(1) \right) Y_w(0)+\nabla_2\phi_i\left(  \bar y(0),\bar y(1) \right) Y_w(1)\leq 0,\quad i\in \hat I_A,
\end{array}
\right.
\end{align*}
and  $w(s)\in T_U^\flat(\bar w(s))$ a.e. $s\in[0,1]$.
 Recalling Definition \ref{sigd},  we see that $(v_0(\cdot), w(\cdot))$ is a singular direction  for problem $(MP)_{[0,1]}$.
 For this $(v_0(\cdot), w(\cdot))$, we use Theorem \ref{sen} to obtain that, there exists a Lagrangian multiplier  $\ell=(\ell_{0_1},\cdots,\ell_{0_r},\ell_1,\cdots,\ell_j,\ell_\psi^\top)^\top$ given by (\ref{lm1}) such that (\ref{snti}) stands, and 
\begin{align}\label{snce}\begin{array}{ll}
&\displaystyle\bar T \int_0^1\Big(  2 \nabla_u H\left(  \bar T s,\bar y(s), p_1^\ell(s), \bar w(s)  \right)(\tilde\sigma(s))
+\nabla_x^2H\left(\bar Ts,\bar y(s), p_1^\ell(s),\bar w(s)  \right)\left(Y_w(s),Y_w(s)\right)
\\&+2\nabla_u\nabla_x H\left(\bar Ts,\bar y(s), p_1^\ell(s),\bar w(s)  \right)\left(Y_w(s), w(s)\right)
+\nabla_u^2H\left(\bar Ts,\bar y(s), p_1^\ell(s),\bar w(s)  \right)\left(w(s),w(s)\right)
\\&-R\left( \tilde p_1^\ell(s), Y_w(s),f(\bar Ts,\bar y(s),\bar w(s)),Y_w(s)  \right)\Big)ds
+D^2\mathcal L_\ell(\bar y(0),\bar y(1))(Y_w)+E(v_0(\cdot),w(\cdot))\leq 0,
\end{array}
\end{align}
holds for all $\tilde\sigma(\cdot)\in\mathcal S_1\stackrel{\triangle}{=}\{\sigma(\cdot)\in L^2(0,1;\mathbb R^m); \sigma(s)\in T_U^{\flat(2)}(\bar w(s), v(s))\;a.e.\,s\in(0,1)\}$, 
where
\begin{align*}
&E(v_0(\cdot),w(\cdot))
\\&= \int_0^1\Big(  \bar T \nabla_t^2H\left(\bar Ts,\bar y(s), p_1^\ell(s),\bar w(s)  \right)y_{v_0}(s)^2
+2\bar T \nabla_{x}\nabla_t H\left(\bar Ts,\bar y(s), p_1^\ell(s),\bar w(s)  \right)\left( y_{v_0}(s), Y_w(s)  \right)
\\&\quad+2y_{v_0}(s)v_0(s)\nabla_t H\left(\bar Ts,\bar y(s), p_1^\ell(s),\bar w(s)  \right)
+2\bar T\nabla_u\nabla_t H\left(\bar Ts,\bar y(s), p_1^\ell(s),\bar w(s)  \right)\left( y_{v_0}(s), w(s) \right)
\\&\quad+2v_0(s)\nabla_x H\left(\bar Ts,\bar y(s), p_1^\ell(s),\bar w(s)  \right)(Y_w(s))
+2v_0(s)\nabla_u H\left(\bar Ts,\bar y(s), p_1^\ell(s),\bar w(s)  \right)(w(s))\Big)ds,
\end{align*}
and $p_1^\ell(\cdot)$ fulfills  (\ref{det}).  Set $p^\ell(t)=p_1^\ell(\frac{t}{\bar T} )$ for $t\in[0,\bar T]$.
By means of (\ref{vtt}),  inequality (\ref{snct1}) holds for all $\sigma(\cdot)\in \mathcal S_{\bar T}$. The proof is completed. $\Box$

\setcounter{equation}{0}
\section{Appendix}\label{ae21}
\def\theequation{4.\arabic{equation}}

To apply Theorems \ref{fon}-\ref{snt} to specific examples, we provide the calculation formulas for some  quantities involved in these theorems within a coordinate chart.

\begin{proposition}\label{ntlc}
Assume $M$ is a complete $n$-dimensional Riemannian manifold with the metric $g$, where $n\in\mathbb N$.  Take a coordinate chart $(\mathcal O, \varphi=(x_1,\cdots,x_n))$ of $M$. For any $x\in \mathcal O$, set $g_{il}(x)=\left\langle\frac{\partial}{\partial x_i},\frac{\partial}{\partial x_l}\right\rangle$ for $i,l=1,\cdots,n$, and let $ G(x) = (g_{il}(x)) $ denote the matrix representation of the metric tensor in this coordinate chart. Then, for a vector $X=\sum_{i=1}^nX_i\frac{\partial }{\partial x_i}$ and a covector $\eta=\sum_{i=1}^n \eta_idx_i$ at $x$, their norms are as follows:
\begin{align}\label{nvc}
|X|=\left|\sqrt{G(x)}\begin{pmatrix}X_1&\cdots&X_n\end{pmatrix}^\top\right|_{\mathbb R^n},\quad |\eta|=\left|\sqrt{ G(x)^{-1}}\begin{pmatrix}
\eta_1&\cdots&\eta_n
\end{pmatrix}^\top\right|_{\mathbb R^n}.
\end{align}
\end{proposition}

{\it Proof.}\;By applying \cite[Exercise 3.8, p.29]{le} to  vector $X$ and covector $\eta$, we can derive (\ref{nvc}).
$\Box$

\begin{proposition}\label{elsn}
Assume that the assumptions in Theorem \ref{sen} hold and  that $(\bar u(\cdot), \bar x(\cdot))$ is a weak Pareto optimal solution of problem $(MP)_{r,T}$ with $\bar u(\cdot)\in L^2(0,T;\mathbb R^m)$. Take a singular direction $v(\cdot)\in L^2(0,T;\mathbb R^m)$, with $X_v(\cdot)$ satisfying (\ref{sdc}) and (\ref{sdc1}). Let $\ell=(\ell_{0_1},\cdots,\ell_{0_r},\ell_1,\cdots,\ell_j,\ell_\psi^\top)^\top$ be a Lagrangian multiplier, with $p^\ell(\cdot)$ satisfying  (\ref{de}). Fix any $t\in[0,T]$. Let $(\mathcal O,\varphi=\xi=(\xi_1,\cdots,\xi_n)^\top)$  be a coordinate chart containing $\bar x(t)$, such that for some $\delta > 0$, we have $\bar x(s) \in \mathcal{O}$ for all $s \in (t - \delta, t + \delta) \cap [0, T]$.  Denote by $\left\{\frac{\partial}{\partial x_i}\Big|_x|\;i=1,\cdots,n\right\}$ and $\left\{dx_i|_x|\;i=1,\cdots,n\right\}$ the bases for the tangent space $T_x M$ and the cotangent space $T_x^*M$, respectively. We also denote by  $\{\Gamma_{pq}^l|\,p,q,l=1,\cdots,n\}$ the set of the  Christoffel symbols corresponding to the Riemannian metric $g$ (see \cite[p.66-67]{pth} for the definition).  For $s\in(t-\delta,t+\delta)\cap[0,T]$, we set
\begin{align*}
&\bar \xi(s)=(\bar \xi_1(s),\cdots,\bar\xi_n(s))^\top=\varphi(\bar x(s)),\quad p^\ell(s)=\sum_{i=1}^n p_i^\ell(s)d\xi_i|_{\bar x(s)},
\\ & \mathbf{p}^{\ell}(s)=(p_{1}^{\ell}(s),\cdots,p_{n}^{\ell}(s))^{\top},\quad X_v(s)=\sum_{i=1}^{n}X_{i}^{v}(s)\frac{\partial}{\partial \xi_{i}}\Big|_{\bar x(s)},
\\ &f(s,\varphi^{-1}(\xi),u)=\sum_{i=1}^n f_i(s,\xi,u)\frac{\partial}{\partial\xi_i}\Big|_{\varphi^{-1}(\xi)},\quad \forall\,(\xi,u)\in\varphi(\mathcal O)\times U,
\\ & H_{\mathcal O}(s,\xi,p_1,\cdots,p_n,u)=\sum_{i=1}^nf_i(s,\xi,u)p_i,\quad \forall\,(\xi,p_1,\cdots,p_n,u)\in\varphi(\mathcal O)\times\mathbb R^n\times U.
\end{align*}
Then, for almost every $s\in(t-\delta, t+\delta)\cap[0,T]$ and $i=1,\cdots,n$, we have
\begin{align}\label{lepxhh}\begin{array}{ll}
&\displaystyle\dot p_i^\ell(s)=-\frac{\partial}{\partial \xi_i}H_{\mathcal O}(s,\bar \xi(s), {\mathbf{p}}^\ell(s),\bar u(s)),
\\[2mm]&\displaystyle \dot X_i^v(s)=\sum_{l=1}^n\frac{\partial}{\partial\xi_l}f_i(s,\bar \xi(s),\bar u(s))X_l^v(s)+\nabla_u f_i(s,\bar \xi(s),\bar u(s))v(s),
\\[2mm]& \displaystyle \nabla_xH(s,\bar x(s), p^\ell(s),w)(X_v(s))=\sum_{i=1}^n\frac{\partial}{\partial\xi_i}H_{\mathcal O}(s,\bar \xi(s),{\mathbf p}^\ell(s),w) X_i^v(s),\quad\forall\,w\in U,
\\[2mm]& \displaystyle \nabla_x^2 H\{s;\ell\}(X_v(s),X_v(s))=\sum_{i,l=1}^n H_{il}(s,\bar\xi(s),{\mathbf p}^\ell(s),\bar u(s))X_i^v(s) X_l^v(s),
\end{array}
\end{align}
where 
\begin{align*}
&H_{il}(s,\bar\xi(s),{\mathbf p}^\ell(s),\bar u(s))
 = \frac{\partial^2}{\partial\xi_i\partial\xi_l} H_{\mathcal O}(s,\bar \xi(s),{\mathbf p}^\ell(s),\bar u(s))-\sum_{\eta=1}^n\Gamma_{li}^\eta(\bar x(s))\frac{\partial}{\partial\xi_\eta}H_{\mathcal O}(s,\bar \xi(s),{\mathbf p}^\ell(s),\bar u(s)).
\end{align*}
Moreover, if $M$ is a two-dimensional manifold,  the following relation holds:
\begin{align}\label{ctt}\begin{array}{ll}
&R(\tilde p^\ell(s), X_v(s),f[s],X_v(s))
=-K(\bar x(s))\left(  p^\ell(s)(f[s])|X_v(s)|^2-p^\ell(s)(X_v(s))\langle f[s],X_v(s)\rangle \right).
\end{array}
\end{align}
Here $K(x)$ is the sectional curvature at $x\in M$ (see \cite[p.83]{pth} for the definition).

\end{proposition}

{\it Proof.}\;
First, we show the first relation of  (\ref{lepxhh}).  By \cite[(3.11)]{dzgJ}, we have
\begin{align}\label{cdpd}
\nabla_{\frac{\partial}{\partial\xi_i}}d\xi_l
=-\sum_{\eta=1}^n\Gamma_{i\eta}^l d\xi_\eta,\;i,l=1,\cdots,n.
\end{align}
Let us fix $s\in(t-\delta, t+\delta)$. Based on the local expression for $\bar x(\cdot)$, we have
 $\dot{\bar x}(s)=\sum_{i=1}^n\dot{\bar \xi}_i(s)\frac{\partial}{\partial\xi_i}\Big|_{\bar x(s)}$. Since $\bar x(\cdot)$ is the solution to (\ref{25}) corresponding to the control $\bar u(\cdot)$, it follows that
\begin{align}\label{ebx}
\dot{\bar \xi}_i(s)=f_i(s,\bar \xi(s),\bar u(s)),\quad i=1,\cdots,n.
\end{align}

For the local expressions for $p^\ell(\cdot)$ and $\dot{\bar x}(\cdot)$, we obtain from   \cite[Proposition 7.2, p.42]{h1}, (\ref{cdpd}),  and (\ref{ebx})   that
\begin{align*}
\nabla_{\dot{\bar x}(s)}p^\ell&=\sum_{i=1}^n \dot p_i^\ell(s)d\xi_i|_{\bar x(s)}-\sum_{i,l,\eta=1}^n p_i^\ell(s) f_l(s,\bar \xi(s),\bar u(s))\Gamma_{l\eta}^i(\bar x(s))d\xi_\eta|_{\bar x(s)},
\end{align*}
which implies
\begin{align}\label{cvdv}
\nabla_{\dot{\bar x}(s)} p^\ell\left(\frac{\partial}{\partial \xi_\eta}\right)=\dot p_\eta^\ell(s)-\sum_{i,l=1}^n p_i^\ell(s)f_l(s,\bar \xi(s),\bar u(s))\Gamma_{l\eta}^i(\bar x(s)),\;\eta=1,\cdots,n.
\end{align}
 Furthermore, we obtain from the definition of the covariant differential (see \cite[p.124]{kn}), the local expressions for $f$ and $p^\ell(\cdot)$, \cite[Proposition 7.2, p.42]{h1} and  the definition for the Christoffel symbols that
\begin{align*}
-\nabla_xf[s]\left(p^\ell(s),\frac{\partial}{\partial\xi_\eta}   \right)
&=-\nabla_{\frac{\partial}{\partial\xi_\eta} }\left(\sum_{i=1}^n f_i(s,\cdot,\bar u(s)) \frac{\partial}{\partial\xi_i} \right)\left( \sum_{l=1}^n p_l^\ell(s)d\xi_l  \right)
\\&=-\frac{\partial}{\partial\xi_\eta}H_{\mathcal O}(s,\bar \xi(s),\textbf{p}^\ell(s),\bar u(s))-\sum_{i,l=1}^n f_i(s,\bar \xi(s),\bar u(s))p_l^\ell(s)\Gamma_{\eta i}^l(\bar x(s)),
\end{align*}
where $\eta=1,\cdots,n$. 
From (\ref{de}), along with the above relation and (\ref{cvdv}), we derive the first formula of (\ref{lepxhh}).

Second, we show the second line of (\ref{lepxhh}). For almost every $s\in(t-\delta,t+\delta)$, we obtain form     \cite[Proposition 7.2, p.42]{h1}, the local expression for $f$ and  the definition of the Christoffel symbols (see \cite[p.66-67]{pth}), that
\begin{align}\label{lecx}
\nabla_{\dot{\bar x}(s)}X_v(d\xi_\eta)= \dot X_\eta^v(s)+\sum_{i,l=1}^n X_i^v(s)f_l(s,\bar \xi(s),\bar u(s))\Gamma_{li}^\eta(\bar x(s)),\;\eta=1,\cdots,n.
\end{align}
We also obtain from the definition of the covariant differential (see \cite[p.124]{kn}) that
\begin{align*}
&\nabla_xf[s](d\xi_\eta, X_v(s))+\nabla_uf[s](d\xi_\eta,v(s))
\\& =\sum_{i=1}^n\frac{\partial}{\partial\xi_i}f_\eta(s,\bar \xi(s),\bar u(s))X_i^v(s)+\sum_{i,l=1}^n X_i^v(s)f_l(s,\bar \xi(s),\bar u(s))\Gamma_{li}^\eta(\bar x(s))
+\nabla_uf_\eta(s,\bar \xi(s),\bar u(s))v(s),
\end{align*}
which together with  (\ref{lecx}), implies the second formula of  (\ref{lepxhh}).

Third,   we show  the last two   formulas of  (\ref{lepxhh}). The third formula follows from \cite[Proposition 2.7, p.123]{kn} and the local expressions for $f(s,\bar x(s),\bar u(s))$ and $p^\ell(s)$.
 We obtain from \cite[Proposition 2.12, p,125]{kn}, the definition of the Christoffel symbol that
\begin{align*}
\nabla_x^2H\{s;\ell\}\left( \frac{\partial}{\partial\xi_i}, \frac{\partial}{\partial\xi_l} \right)
&=\frac{\partial}{\partial\xi_l}\left( \frac{\partial}{\partial\xi_i} (H)   \right)-\nabla_{ \frac{\partial}{\partial\xi_l}} \frac{\partial}{\partial\xi_i}(H),\;i,l=1,\cdots,n,
\end{align*}
which implies the last line of  (\ref{lepxhh}). 

To prove (\ref{ctt}), consider the curvature tensor $\tilde R$ as defined in \cite[p.78]{pth}: $
\tilde R(X,Y)Z=\nabla_X\nabla_Y Z-\nabla_Y\nabla_XZ-\nabla_{[X,Y]}Z$, 
over vector fields $X,Y,Z$, with $[X,Y]\stackrel{\triangle}{=}YX-XY$. The tensor field $\tilde R$ is a  $(1,3)$-tensor field. Using the Riemannian metric $g$, the tensor field $\tilde R$ can be changed to a $(0,4)-$tensor field as follows:
$
R(X,Y,Z,W)\stackrel{\triangle}{=}\langle \tilde R(X,Y)Z,W\rangle$,
where $X,Y,Z,W$ are vector fields. Then tensor field $R$ is also called the curvature tensor (see \cite[p.79]{pth}). At the same time, we denote that
$
(X\wedge Y)Z=\langle X,Z\rangle Y-\langle Y, Z\rangle X,
$
for vector fields $X,Y,Z$ (see \cite[p.82]{pth}).
Then, according to \cite[Proposition 3.1.3, p.84]{pth}, we can obtain
\begin{align*}
R\left( \tilde p^\ell(s), X_v(s), f[s], X_v(s)  \right)
&=-K(\bar x(s))\left\langle  (\tilde p^\ell(s)\wedge X_v(s))f[s],X_v(s)  \right\rangle
\\&=-K(\bar x(s))\left\langle \langle\tilde p^\ell(s),f[s]\rangle X_v(s)-\langle X_v(s),f[s]\rangle \tilde p^\ell(s),X_v(s)  \right\rangle,
\end{align*}
which implies (\ref{ctt}).
The proof is completed. $\Box$

Next, we consider a special class of Riemannian manifolds, which is the graph of
a function.
\begin{proposition}\label{csg1}
Let $a:\mathbb R^2\to\mathbb R$ be a $C^\infty$ map. Then, the graph of function $a$, given by $M=\{(x_1,x_2,a(x_1,x_2))^\top|\,  (x_1,x_2)\in\mathbb R^2\}$,  is  a two-dimensional manifold with the coordinate chart
\begin{align}\label{ccg}
\varphi:(x_1,x_2,a(x_1,x_2))^\top(\in M)\mapsto (x_1,x_2)^\top.
\end{align}
In this chart, the basis for the tangent space of $M$ at $(x_1,x_2,a(x_1,x_2))\in M$ consists of
\begin{align}\label{btg}
\frac{\partial}{\partial x_1}\Big|_M=\begin{pmatrix}
1&0&a_{x_1}
\end{pmatrix}^\top,\quad \frac{\partial}{\partial x_2}\Big|_M=\begin{pmatrix}
1&0&a_{x_2}
\end{pmatrix}^\top,
\end{align}
where
$a_{x_i}$ denotes the  partial derivative of function $a$ with respect to the variable $x_i$.
Furthermore, we endow a Riemannian metric $g$ (see \cite[Definition 2.1, p.38]{c} for the definition), which can be represented  in this 
chart by $g_{il}=\frac{\partial}{\partial x_i}\Big|_M \cdot \frac{\partial}{\partial x_l}\Big|_M$  ($i,l=1,2$), where the dot ``$\cdot$'' represents the dot product of two vectors in Euclidean space 
$\mathbb R^3$.  These components are given by
\begin{align}\label{rmg}
\begin{array}{ll}
& \displaystyle g_{11}=1+a_{x_1}^2,\quad g_{12}=g_{21}=a_{x_1}a_{x_2},
\quad g_{22}=1+a_{x_2}^2.
\end{array}
\end{align}
 Then, the Christoffel symbols are as follows:
\begin{align}\label{csg}\begin{array}{l}
\displaystyle\Gamma_{11}^1=\frac{a_{x_1x_1}a_{x_1}}{1+a_{x_1}^2+a_{x_2}^2},\quad \Gamma_{11}^2=\frac{a_{x_1x_1}a_{x_2}}{1+a_{x_1}^2+a_{x_2}^2},\quad
 \Gamma_{12}^1=\Gamma_{21}^1=\frac{a_{x_1}a_{x_1x_2}}{1+a_{x_1}^2+a_{x_2}^2},
 \\[2mm]\displaystyle \Gamma_{12}^2=\Gamma_{21}^2=\frac{a_{x_2}a_{x_1x_2}}{1+a_{x_1}^2+a_{x_2}^2},\quad
\Gamma_{22}^1=\frac{a_{x_2x_2}a_{x_1}}{1+a_{x_1}^2+a_{x_2}^2},\quad \Gamma_{22}^2=\frac{a_{x_2x_2}a_{x_2}}{1+a_{x_1}^2+a_{x_2}^2},
\end{array} 
\end{align}
and the sectional curvature at $x\in M$ is given by
\begin{align}\label{scg}
K(x)=\frac{a_{x_1x_1}a_{x_2x_2}-a_{x_1x_2}^2}{(1+a_{x_1}^2+a_{x_2}^2)^2}.
\end{align}
\end{proposition}

{\it Proof.}\;  It follows from \cite[Proposition 1, p.58]{c1} and its proof  that $M$ is a two-dimensional manifold, equipped with the coordinate chart $\varphi$. Furthermore, according to  \cite[p.84]{c1}, we find that  $\left\{\frac{\partial}{\partial x_1}\Big|_M, \frac{\partial}{\partial x_2}\Big|_M \right\}$
is a basis for the tangent space at each point of $M$.

Using the formula for Christoffel symbols (\cite[p.56]{c}), a direct calculation yields (\ref{csg}).
According to \cite[Remark 2.7, p.131]{c}, 
 for a two-dimensional Riemannian manifold, the sectional curvature is identical to the Gaussian curvature. Therefore, from \cite[Example 5, p.162]{c1}, we derive the formula for the sectional curvature given in equation (\ref{scg}). This completes the proof. $\Box$

Finally, we prove Proposition \ref{psd}.

{\it Proof of Proposition \ref{psd}.}\; Let $\ell$ be a Lagrangian multiplier, and let  $p^\ell(\cdot)$ be the solution to (\ref{de}) corresponding to the vector $\ell$.
 We obtain from  (\ref{sdc}), (\ref{sdc1}), (\ref{de}), \cite[(5.9)]{dzgJ} and integration by parts that,
\begin{align*}
0=&\nabla_2\mathcal L_{\ell}(\bar x(0),\bar x(T))(X_v(T))+\nabla_1\mathcal L_{\ell}(\bar x(0),\bar x(T))(X_v(0))
\\=& p^\ell(T)(X_v(T))-p^\ell(0)(X_v(0))
\\=&\int_0^T\left( \nabla_{\dot{\bar x}(t)}  p^\ell(X_v(t))+p^\ell(t)(\nabla_{\dot{\bar x}(t)}X_v)\right) dt
\\=&\int_0^T\nabla_u H(t,\bar x(t), p^\ell(t),\bar u(t))(v(t)) dt,
\end{align*}
which, together with (\ref{lm3}), completes the proof. $\Box$

\medskip

\bibliographystyle{amsplain}

\end{document}